\documentclass[preprint,1p,times]{elsarticle}

\usepackage{graphicx}

\usepackage{amssymb}
\usepackage{amsmath}

\usepackage{float}

\usepackage{color}
\definecolor{red}{rgb}{1,0,0}
\definecolor{green}{rgb}{0,0.5,0}
\definecolor{blue}{rgb}{0,0,1}

\journal{Journal of Computational Physics}

\newcount\ndots
\def\drawline#1#2{\raise 2.5pt\vbox{\hrule width #1pt height #2pt}}
\def\spacce#1{\hskip #1pt}

\def\chndash{\hbox {\drawline{8.5}{.5}\spacce{2}\drawline{3}{.5}\spacce{2}\drawline{8.5}{.5}}\nobreak\ }

\def\chndash2{\hbox {\drawline{4}{.5}\spacce{2}\drawline{1.5}{.5}\spacce{2}\drawline{4}{.5}}\nobreak\ }
\def\trian{\raise 1.25pt\hbox{$\scriptscriptstyle\triangle$}\nobreak\ }
\def\circle{$\circ$\nobreak\ }

\def\square{${\vcenter{\hrule height .4pt
        \hbox{\vrule width .4pt height 3pt \kern 3pt
        \vrule width .4pt}
        \hrule height .4pt}}$\nobreak\ }
\def\plus{\raise 1.25pt \hbox{$\scriptscriptstyle +$}\nobreak\ }
\def\x{\raise 1.25pt \hbox{$\scriptscriptstyle \times$}\nobreak\ }

\def\solidtrian{\raise 1.25pt
   \hbox to 3bp{
\hfill}\nobreak\ }
\def\solidsquare{\vrule height .9ex width .8ex depth -.1ex\nobreak\ }

\def\solidcclose{\drawline{10}{.5}\nobreak\raise
  0.5pt\hbox{$\bullet$}\drawline{10}{.5}\nobreak\ }

\def\solidsclose{\drawline{10}{.5}\nobreak\raise
  0.5pt\hbox{\solidsquare}\drawline{10}{.5}\nobreak\ }

\def\solidtclose{\drawline{10}{.5}\nobreak\raise
  0.5pt\hbox{\solidtrian}\drawline{10}{.5}\nobreak\ }

\def\solidcopen{\drawline{10}{.5}\nobreak\raise
  0.5pt\hbox{\circle}\drawline{10}{.5}\nobreak\ }

\def\solidsopen{\drawline{10}{.5}\nobreak\raise
  0.5pt\hbox{\square}\drawline{10}{.5}\nobreak\ }

\def\solidtopen{\drawline{10}{.5}\nobreak\raise
  0.5pt\hbox{\trian}\drawline{10}{.5}\nobreak\ }

\def\solidx{\drawline{10}{.5}\nobreak\raise
  0.5pt\hbox{\x}\drawline{10}{.5}\nobreak\ }

 \def\be{\begin{eqnarray}} 
 \def\ee{\end{eqnarray}}  
 \def\bfi{\begin{figure}[htb]} 
 \def\efi{\end{figure}} 
  \def\bta{\begin{table}[htb] \small}  
  \def\enta{\end{table}}

\def\ro0{\rho_0} 
\def\x{{\bf x}} 
\def\rh1{\rho_1}  
   
\usepackage{subfigure}
\usepackage{subfigmat}
\usepackage{diagbox}

\usepackage{array}
\newcolumntype{L}[1]{>{\raggedright\let\newline\\\arraybackslash\hspace{0pt}}m{#1}}
\newcolumntype{C}[1]{>{\centering\let\newline\\\arraybackslash\hspace{0pt}}m{#1}}
\newcolumntype{R}[1]{>{\raggedleft\let\newline\\\arraybackslash\hspace{0pt}}m{#1}}

\begin{document}

\begin{frontmatter}



\title{A Hybrid Adaptive Low-Mach-Number/Compressible Method: Euler Equations.}


\author[CCSE]{Emmanuel Motheau\corref{cor1}}
\ead{emotheau@lbl.gov}
\author[CCSE]{Max Duarte}
\author[CCSE]{Ann Almgren}
\author[CCSE]{John B.  Bell}

\address[CCSE]{Center for Computational Sciences and Engineering, Computational Research Division, Lawrence Berkeley National Laboratory, Berkeley, CA 94720-8139, USA.}

\cortext[cor1]{Corresponding author.}

\begin{abstract}
Flows in which the primary features of interest do not rely on high-frequency acoustic effects, but in which long-wavelength acoustics play a nontrivial role, present a computational challenge. Integrating the entire domain with low-Mach-number methods would remove all acoustic wave propagation, while integrating the entire domain with the fully compressible equations can in some cases be prohibitively expensive due to the CFL time step constraint. For example, 
simulation of thermoacoustic instabilities might require fine resolution of the fluid/chemistry interaction but not require fine resolution of acoustic effects, yet one does not want to neglect the long-wavelength wave propagation and its interaction with the larger domain.

The present paper introduces a new multi-level hybrid algorithm to address these types of phenomena. In this new approach, 
the fully compressible Euler equations are solved on the entire domain, potentially with local refinement, 
while their low-Mach-number counterparts are solved on subregions of the domain with higher spatial resolution. 
The finest of the compressible levels communicates inhomogeneous divergence constraints to the coarsest of the low-Mach-number levels, 
allowing the low-Mach-number levels to retain the long-wavelength acoustics.  The performance of the hybrid method is shown for a series of test cases, including results from a simulation of the aeroacoustic propagation generated from a Kelvin-Helmholtz instability 
in low-Mach-number mixing layers. It is demonstrated that compared to a purely compressible approach, 
the hybrid method allows time-steps two orders of magnitude larger at the finest level, leading to an overall reduction of the computational time by a factor of $8$.
\end{abstract}


\end{frontmatter}



\section{Introduction}

Many interesting fluid phenomena occur in a regime in which the fluid velocity is much less than the speed of sound. Indeed, it is possible to make a distinction between scales of fluctuations, depending on how a hydrodynamic fluid element is sensitive to acoustic perturbations. Acoustic waves that do not carry enough energy to perturb a flow are referred to short-wavelengths. In contrary, long-wavelengths refer to large scale motions where acoustic and hydrodynamic fluctuations can interact. Low-Mach-number \cite{Majda:1985,Giovangigli:1999,Knikker:2011} schemes exploit the separation of scales between acoustic and advective motions; these methods calculate the convective flow field but do not allow explicit propagation of acoustic waves. Their computational efficiency relative to explicit compressible schemes results from the fact that the time step depends on the fluid velocity rather than sound speed. However, there is a class of problems for which the small-scale motions can be adequately captured with a low-Mach-number approach, but which require in addition the representation of long wavelength acoustic waves.
This paper introduces a computational methodology for accurately and efficiently
calculating these flows.

An important example of this type of flow is thermoacoustic instabilities in 
large scale gas turbine engines. In these engines the region where the burning takes place can be modeled using a low-Mach-number approach, since the short-wavelength acoustic waves generated by the heat release do not carry sufficient information or energy to be of interest. Low-Mach-number modeling of turbulent combustion has been demonstrated to be an efficient way to generate accurate solutions \cite{McMurtry:1986,Day:2000,Najm:2010,Aspden:2016,Motheau:2016a,Wang:2017}.  However, in large burners, under certain conditions the long-wavelength acoustic waves that emanate from the burning region can reflect from the walls of the burner and impinge on the burning region, generating thermoacoustic instabilities which can be violent enough to disrupt the flame, as well as lead to mechanical failures or excessive acoustic noise \cite{Ducruix:2003p739,Huang:2009,Gicquel:2012,Lieuwen:2012book,Candel:2013,Motheau:2014b}. There is currently a great deal 
of interest in the problem of how to control the instabilities through
passive or active control mechanisms \cite{Poinsot:2017}.

This scenario could clearly be modeled using the fully compressible reacting flow equations, but the sound speed is high and the burners are large, and performing such a simulation at the resolution required for detailed characterization of the flame is computationally infeasible. Thus the goal of the work here is to construct a methodology in which the time scale at which the equations are evolved is that of the fluid velocity rather than the sound speed, but which can explicitly propagate the long-wavelength acoustic waves as they travel away from the flame and as they return and interact with the flame that created them.

This paper is the first of a series of papers describing the development
of this methodology. For the purposes of this paper, one of the simplest 
low-Mach-number equation sets is considered, i.e. the variable density incompressible Euler equations. These equations allow different regions of the flow to have different densities, but do not allow any volumetric changes to occur (i.e. the material derivative of the density is zero). A hybrid approach is constructed in which variants of both the low-Mach-number equations and the fully compressible equations are solved in each time step; the computational efficiency of this approach results from the fact that the compressible equations are solved at a coarser resolution than the low-Mach-number equations. As a result, only long
wavelength acoustic waves are resolved, yet the fine scale locally incompressible structure can still be resolved on the finer level(s).   

The method is similar to the Multiple Pressure Variables (MPV) first introduced in a set of papers by Munz \emph{et al} \cite{Roller:2000,Munz:2003,Park:2005,Munz:2007}. The essence of the MPV approach is to decompose the pressure into three terms: the thermodynamic pressure $p_0$; the acoustic pressure $p_1$; and the perturbational pressure $p_2$. The acoustic signal is carried by $p_1$, and $p_2$ is used to satisfy the divergence constraint on the low-Mach-number levels and is defined as the solution to a Poisson equation. Different approaches for solving $p_1$ were proposed in the aforementioned references, for example by solving a set of Linearized Euler Equations (LEEs) on a grid that is a factor of $1/M$ coarser, where $M$ is a measure of the Mach number of the flow.  
Differently, \citet{PeetLele:2008} developed a hybrid method in which the exchange of information between the fully compressible and 
low-Mach-number regions occurs through the boundary conditions of overlapping meshes. 
The novelty of the present paper is that the fully compressible equations are solved without any approximation, 
and that an adaptive mesh refinement (AMR) framework is employed to optimize the performance of the algorithm. 
Thus, while the fully compressible equations are solved in the entire domain, with possible additional local refinement, 
the hybrid strategy developed in the present paper allows refined \emph{patches} where the low-Mach-number equations 
are solved at finer resolution.

Note that there have been a number of other approaches to bridging the gap between fully compressible and low-Mach-number approaches. One alternative to the MPV methodology are the so-called \emph{unified}, \emph{all-speed}, \emph{all-Mach} or \emph{Mach-uniform} approaches \cite{Yoon:1999,vanderHeul:2003,Nerinckx:2005,Cordier:2012}, which consist of a single equation set that is valid from low to high Mach numbers. These methods retain the full compressible equation set, but numerically separate terms which represent convection at the fluid speed from acoustic effects traveling at the sound speed. Inherent in these approaches is that at least some part of the acoustic signal is solved for implicitly, which makes them inapplicable for our applications of interest in which explicit propagation of the long wavelength acoustic modes is preferred. 

Note also that all of the methods described above involve feedback from the compressible solution to the low-Mach-number solution, and the reverse, thus they fundamentally differ from many hybrid methods employed in the aeroacoustics community, in which the acoustic calculation does not feed back into the 
low-Mach-number solution. Methods such as Expansion about Incompressible Flow (EIF) 
\cite{HardinPope:1994} can be used to calculate acoustic waves via Lighthill's analogy approach given an existing incompressible solution. A review of aeroacoustic methods is beyond the scope of this paper, but a comparison of EIF, MPV and LEEs is given in \citet{RollerEtAl:2005}. More recently, many groups \cite{Bogey:2002,Fortune:2004,Golanski:2004,Golanski:2005} have investigated the coupling between a low-Mach-number detailed simulation of noise sources from a small scale turbulent flow, and the aeroacoustic propagation within a larger domain with the LEEs. It will be shown in the results section that the novel hybrid method developed in the present paper is able to tackle the same kind of problem 
while solving the purely compressible equations instead of the LEEs and allowing feedback of the acoustics into the low-Mach-number solution.

The remainder of this paper is organized as follows.  In Section~\ref{sec:equations}
the hybrid hierarchical grid strategy and governing equations that are
solved at each resolution are presented. Then, in Section~\ref{sec:methods_num} the
time advancement algorithm is detailed, as well as the procedures for interpolation and exchange of the variables between the different sets of equations 
at different levels.  
Finally, Section~\ref{sec:results} contains the numerical results of the canonical test cases employed to assess 
the spatial and temporal rates of convergence of the hybrid method, 
as well as the simulation of the propagation of aeroacoustic waves generated by the formation of a Kelvin-Helmholtz instability in mixing layers.  Note that these numerical examples are computed in 2D, but it is emphasized that the algorithm presented in this paper can be easily extended to 3D.

\section{Hybrid hierarchical grid strategy and governing equations}
\label{sec:equations}

The key idea of the algorithm developed in the present paper is to separate the acoustic part of the flow from the hydrodynamics, and to retain acoustic effects only at wavelengths at longer length scales than the finest flow features.  This is achieved by solving a modified form of the low-Mach-number equations at the resolution required by the fine scale features of the flow, while solving the fully compressible governing equations on a coarser level (or levels) underlying the low-Mach-number levels. Because the compressible equations are not solved at the finest level, the overall time step is reduced by a factor of the ratio of grid resolutions from what it would be in a uniformly fine compressible simulation. It is important to note here that $\Delta t_{\rm{LM}}/\Delta t_{\rm Comp} \approx 1/M$, where $\Delta t_{\rm Comp}$ and $\Delta t_{\rm{LM}}$ are the time-steps associated to the fully compressible and low-Mach-number equations. If a ratio of $2$ in resolution is considered between the compressible and low-Mach-number levels, this means that the advancement of the fully compressible equations will be performed with a number of sub-steps scaling with $1/(2M)$. Consequently, $\Delta t_{\rm Comp}$ and $\Delta t_{\rm{LM}}$ will be virtually the same for Mach numbers $M > 0.5$. In other terms, the numerical strategy developed in the present paper is not suitable to be applied in regions of flows featuring a Mach number above a value of $0.5$. Moreover, for Mach numbers in the range of $0.25<M<0.5$, one iteration performed over the low-Mach-number level would involve the time advancement of the compressible equations within two time-steps on the coarser level. As the present algorithm involves a projection method with successive solve of a Poisson equation, the additional computational cost may not be interesting compared to the advancement of the equations with a purely compressible method. Consequently, in practice, it is estimated that the present numerical strategy is valuable and represents a gain in computational time when applied in regions of flows that feature Mach numbers $M<0.2$.

In practice, the grid hierarchy can contain multiple levels for each of the two solution approaches. This fits naturally within the paradigm of block-structured adaptive mesh refiment (AMR), although most published examples of AMR simulations solve the same set of equations at every level.  The present algorithm forms the \textbf{LAMBDA} code and is developed within the BoxLib package \cite{BoxLibref:2016,BoxLib}, a hybrid C++ /Fortran90 software framework that provides support for the development of parallel structured-grid AMR applications. 

The computational domain is discretized into one or more grids on a set of different levels of resolution.  The levels are denoted by
$l=1, \cdots, L$. The entire computational domain is covered by the coarsest level ($l=1$); the finest level is denoted by $l = L.$
The finer levels may or may not cover the entire domain; the grids at each level are properly nested in the sense that the 
union of grids at level $l+1$ is contained in the union of grids at level $l$. 
The fully compressible equations are solved on the \emph{compressible levels}, which are denoted as 
$l_{\rm Comp} = \left\{1,\cdots,l_{\rm max\_comp}\right\}$, while on the \emph{low-Mach levels} denoted as 
$l_{\rm LM} = \left\{l_{{\rm max\_comp}+1},\cdots,L\right\}$, 
the modified low-Mach-number equations are solved. The index ${\rm  max\_comp}$ is an integer that denotes here the total number of compressible layers involved in the computation. 
For ease of implementation of the interpolation procedures, the current algorithm assumes a ratio of $2$ in resolution between 
adjacent levels and that the cell size on each level is independent of direction.

\subsection{Governing equations solved on compressible level}
\label{subsec:eqs_on_compressible_layers}

The set of fully compressible Euler equations are solved on levels $l_{\rm Comp} = \left\{1,\cdots,l_{\rm max\_comp}\right\}.$
The conservation equations for continuity, momentum and energy are expressed as:
\begin{align}
  \frac{\partial \rho }{\partial t} &=  - \nabla \cdot \left(\rho \mathbf{u} \right) 
\label{eqn:mass_comp} \\
  \frac{\partial \left(\rho \mathbf{u} \right) }{\partial t} &= - \nabla \cdot \left(\rho \mathbf{u} \mathbf{u} 
\right) - \nabla p_{\rm Comp} \label{eqn:momentum_comp} \\
  \frac{\partial \left(\rho E\right) }{\partial t} &= - \nabla \cdot \left(\rho
\mathbf{u} E + p_{\rm Comp} \mathbf{u} \right) \label{eqn:energy_comp}
\end{align}

Here, $\rho$, $\mathbf{u}$ and $E$ are the mass density, the velocity vector and the total energy per unit mass, respectively. The total energy is expressed as $E = e + \mathbf{u} \cdot 
\mathbf{u} / 2$, where $e$ is the specific internal energy. The total pressure $p_{\rm Comp}$ is 
related to the energy through the following equation of state:
\begin{equation}
p_{\rm Comp} = \left( \gamma - 1 \right) \rho e
\label{eqn:eos}
\end{equation}
where $\gamma$ is the ratio of the specific heats. Note that Eq.~(\ref{eqn:eos}) represents a very simplified assumption of the equation of state, and that it will be generalized in future work, for example to deal with reactive Navier-Stokes equations composed of multiple chemical species.

\subsection{Governing equations solved on low-Mach levels}
\label{subsec:eqs_on_LM_layers}

The set of governing equations are recast under the low-Mach-number assumption 
and solved on levels $l_{\rm LM} = \left\{l_{{\rm max\_comp}+1},\cdots,L\right\}$. The description of the mathematical derivation of the equations under this assumption is out of scope of the present paper, and can be found in the seminal works of \citet{Majda:1985} and \citet{Giovangigli:1999}. However, from a numerical point of view, it should be noted that different ways to arrange the conservation equations are possible, but as recalled by \citet{Knikker:2011} in his review paper, it is not possible to solve all of them in a conservative form unless an implicit approach is employed. As it will be detailed in \S\ref{subsubsec:step5}, the present algorithm is based on the strategy initially proposed by \citet{Day:2000}, which aims to advance the mass and energy equations while enforcing the conservation of the equation of state through a modification of the constraint on the divergence. In summary, here in the present algorithm mass and energy are formally conserved, while the momentum is conserved up to $\mathcal{O}\left(2\right)$ accuracy. The conservation equations for continuity, momentum, and energy are, respectively:

\begin{align}
  \frac{\partial \rho }{\partial t} &=  - \nabla \cdot \left(\rho \mathbf{u} \right) \label{eqn:mass_LM} \\
  \frac{\partial \mathbf{u} }{\partial t} &=  - \mathbf{u} \cdot \nabla \mathbf{u} - 
\frac{1}{\rho} \nabla \left(p_0 + p_1 + p_2 \right) \label{eqn:momentum_LM} \\
\frac{\partial \left( \rho h \right) }{\partial t} &=  - \nabla \cdot \left(\rho 
\mathbf{u} h  \right) + \frac{{\rm D} p_1}{{\rm D} t} \label{eqn:energy_LM}
\end{align}
where $h=e+p/\rho$ is the enthalpy. Eqs.~(\ref{eqn:mass_LM})-(\ref{eqn:energy_LM}) are accompanied by the following constraint on the velocity:
\begin{equation}
\nabla \cdot \mathbf{u} = \nabla \cdot \mathbf{u}_{\rm Comp}
\end{equation}
where $\mathbf{u}_{\rm Comp}$ is 
interpolated from the compressible level. As explained in the introduction, the pressure that appears in the low-Mach-number equations 
is not written as a single term like $p_{\rm Comp}$ in the fully compressible equations, but has been decomposed into three terms: the thermodynamic pressure $p_0$, the acoustic pressure $p_1$, and the perturbational pressure $p_2$. As will be explained below in the full description of the integration algorithm, $p_0$ is constant through the whole simulation, while $p_1$ is provided from the compressible solution and $p_2$ is intrinsic to the projection method for the pressure. It should be noted that these pressure terms are derived quantities from the mass and the enthalpy, which are conserved quantities advanced in time with Eqs.~(\ref{eqn:mass_LM}) and~(\ref{eqn:energy_LM}). Thus, one should emphasize that the density is not decomposed during the projection procedure. Following on, the pressure terms described above are derived quantities from the mass and the enthalpy. In the standard low-Mach-number approximation it is the background pressure $p_0$ that satisfies the equation of state. In the current model in which the low-Mach-number equations incorporate long wavelength acoustics, it is the sum of the background $p_0$ and hydrodynamic pressure $p_1$ that satisfy the equation of state; see Eq.~(\ref{eqn:init_energy_LM}). The mathematical description of the algorithm for the time integration is presented below.

\section{Integration procedure}
\label{sec:methods_num}

\subsection{Overall presentation of the algorithm}
\label{subsec:overall_presentation}

At the beginning of a time-step, both the compressible and the low-Mach-number equations share the same state variables on all levels. 
The procedure can be summarized as follows:

\begin{enumerate}
  \item The time-steps for the fully compressible equations as well as the low-Mach-number equations have to be computed and synchronized first so as to define a global time-marching procedure. 
  \item The fully compressible Eqs.~(\ref{eqn:mass_comp}-\ref{eqn:energy_comp}) are advanced in time on the designated compressible levels through the whole time-step, from $t^n$ to $t^{n+1}$. As explained at \S\ref{subsubsec:step1}, this may involve several sub-steps depending on the flow and mesh configurations. At the end of the procedure, state variables are known on those levels at $t^{n+1}$.
  \item The low-Mach-number Eqs.~(\ref{eqn:mass_LM}-\ref{eqn:energy_LM}) are then advanced in time on the designated low-Mach levels from $t^n$ to $t^{n+1}$. The terms involving the acoustic pressure $p_1$ are provided by interpolation from the compressible solution. As the momentum Eq.~(\ref{eqn:momentum_LM}) is advanced through a fractional-step method, a variable-coefficient Poisson equation must be solved to correct the velocity fields. 
The constraint on the velocity that appears as a source term in the Poisson equation is provided by construction with interpolated values from the compressible solution. At the end of the procedure, state variables on the low-Mach levels are spatially averaged down to the compressible levels and a new time-step can begin.
\end{enumerate}

The algorithm detailed below constitutes the new \textbf{LAMBDA} code, and uses routines from the existing codes \textbf{CASTRO} \cite{Almgren:2010a} 
and \textbf{MAESTRO} \cite{Maestro_ref:2010}. This ease of reuse and demonstrated accuracy of the existing discretizations 
motivated the choices of the numerical methods described in the present paper; 
however, the algorithm presented here could be adapted to use alternate discretizations.

\subsection{Temporal integration}
\label{subsec:methods_num_algo_global}

At the beginning of a simulation, the density $\rho^{\rm init}$, the velocity vector $\mathbf{u}^{\rm init}$ and total pressure $p_{\rm Comp}^{\rm init}$ 
are specified as the initial conditions.  The pressure $p_{\rm Comp}^{\rm init}$ is specified as the sum of a static reference pressure $p_0^{\rm init}$, 
which will remain constant through the whole simulation, and a possible acoustic fluctuation $p_1^{\rm init}$ that depends on the initial solution. 

The variables on the {\em compressible} levels are initialized as 
\begin{align}
\rho &= \rho^{\rm init} \label{eqn:init_density_comp} \\
\rho \mathbf{u} &= \rho^{\rm init}\mathbf{u}^{\rm init} \label{eqn:init_momentum_comp} \\
\rho E &= \frac{p_0^{\rm init}+p_1^{\rm init}}{\gamma -1} + \frac{1}{2}\rho^{\rm init} \mathbf{u}^{\rm init}\cdot\mathbf{u}^{\rm init} \label{eqn:init_energy_comp}
\end{align}
and those on the {\em low-Mach-number} levels are initialized as
\begin{align}
      \rho &= \rho^{\rm init} \label{eqn:init_density_LM} \\
\mathbf{u} &= \mathbf{u}^{\rm init} \label{eqn:init_momentum_LM} \\
    \rho h &=  \left(p_0^{\rm init}+p_1^{\rm init}\right) \left(1 + \frac{1}{\gamma -1} \right) \label{eqn:init_energy_LM}
\end{align}
 
\subsubsection{Step 1: Computation of time-steps}
\label{subsubsec:step1}

The very first step of the time-integration loop is to compute the time-steps $\Delta t_{\rm Comp}$ and $\Delta t_{\rm{LM}}$ associated to the fully compressible and low-Mach-number equations, respectively:

\begin{align}
  \Delta t_{\rm Comp} &= \sigma^{{\rm CFL}} \min_{l_{\rm Comp}} \left\{ \frac{\Delta x}{|\mathbf{u}| + c} \right\} \\
  \Delta t_{\rm{LM}} &= \sigma^{{\rm CFL}} \min_{l_{\rm LM}} \left\{ \frac{\Delta x}{|\mathbf{u}| } \right\}
\end{align}
where $\min_{l_{\rm Comp}}$ and $\min_{l_{\rm LM}}$ are the minimum values taken over all computational grid cells that belong to the set of levels $l_{\rm Comp}=\left\{1,\cdots,l_{\rm max\_comp}\right\}$ and $l_{\rm LM}=\left\{l_{{\rm max\_comp}+1},\cdots,L\right\}$, respectively. The CFL condition number $0 < \sigma^{{\rm CFL}} < 1$ is set by the user, and $c=\sqrt{\gamma p_{\rm Comp} / \rho}$ is the sound speed computed with the pressure coming from the fully compressible equations. Note here that for the ease of implementation and presentation, the algorithm does not employ the specific AMR technique of sub-cycling in time between levels where the same equations are solved. It is emphasized that the hybrid strategy can be easily adapted to such technique.

The particularity of the present hybrid algorithm is that the resolution of the low-Mach-number level(s) is always finer than the finest compressible level. However, the time-step for evolving the low-Mach-number equations depends on the flow velocity, while the compressible time-step depends on both the flow velocity and the sound speed. Thus, one has to guarantee that the low-Mach-number time-step is not smaller than the compressible time-step, viz. $\Delta t_{\rm Comp} \leqslant \Delta t_{\rm{LM}}$. Consequently, depending on the local sound speed, the time-advancement of the fully compressible equations may involve several sub-steps $K$, and an effective hybrid time-step is defined as:
\begin{equation}
  \Delta t_{\rm hyb} = \frac{\Delta t_{\rm{LM}}}{K} \label{eqn:delta_t_hyb}
\end{equation}
with 
\begin{equation}
  K = \Bigg\lceil \frac{\Delta t_{\rm{LM}}}{\min\left(\Delta t_{\rm Comp},\Delta t_{\rm{LM}} \right)} \Bigg\rceil
  \label{eqn:ceiling_time_steps}
\end{equation}
Note that in Eq.~(\ref{eqn:ceiling_time_steps}), $\lceil \cdot \rceil$ is the ceiling function.

\subsubsection{Step 2: Time advancement of the fully compressible equations}
\label{subsubsec:step2}

Recall that 
the fully compressible conservative Eqs.~(\ref{eqn:mass_comp}-\ref{eqn:energy_comp}) are advanced in time from 
$t^n$ to $t^{n+1}$ through $K$ sub-steps of $\Delta t_{\rm hyb}$, and for all levels $l_{\rm Comp}=\left\{1,\cdots,l_{\rm max\_comp}\right\}$. The integration procedure during this step is complex and will only be summarized below. Note that as the present \textbf{LAMBDA} code is directly reusing routines from the \textbf{CASTRO} code \cite{Almgren:2010a} for the integration of the fully compressible equations, the algorithm is summarized below and the reader is referred to the {\bf CASTRO} references for additional detail.

Eqs.~(\ref{eqn:mass_comp}-\ref{eqn:energy_comp}) are solved in their conservative form as follows:

\begin{equation}
  \mathbf{U}^{k+1} = \mathbf{U}^{k} - \Delta t_{\rm hyb} \nabla \cdot \mathbf{F}^{k+1/2} 
\label{eqn:compressible_vectorised}
\end{equation}
where $k=0,\ldots,K-1$. Here $\mathbf{U}$ is the conserved state vector (stored at cell-centers) and  $\mathbf{F}$ is the flux vector (located at edges of a cell):

\begin{equation}
  \mathbf{U}=\left\{ \begin{array}{c}
                      \rho \\
		      \rho \mathbf{u} \\
		      \rho E
                     \end{array} \right\}
 \label{eqn:vector_conservative_vars}
\end{equation}
and  
\begin{equation}
  \mathbf{F}=\left\{ \begin{array}{c}
                      \rho \mathbf{u} \\
		      \rho \mathbf{u}\mathbf{u}-p_{\rm Comp} \\
		      \rho \mathbf{u} E - p_{\rm Comp}\mathbf{u}
                     \end{array} \right\}
\end{equation}
Note that at the beginning of the first sub-step, $\mathbf{U}^{k=0}=\mathbf{U}^{n}$. Similarly, at the end of the last sub-step, $\mathbf{U}^{n+1}=\mathbf{U}^{K}$.

The edge-centered flux vector $\mathbf{F}^{k+1/2}$ is constructed from time-centered edge states computed with a conservative, shock-capturing, unsplit Godunov method, which makes use of the
Piecewise Parabolic Method (PPM), characteristic tracing and full corner coupling \cite{Miller:2002,Colella:2008,Almgren:2010a}. Basically this particular procedure follows four major steps:

\begin{enumerate}
 \item The conservative Eqs.~(\ref{eqn:mass_comp}-\ref{eqn:energy_comp}) are rewritten in terms of the primitive state vector,
$\mathbf{Q}=\left\{\rho, \mathbf{u}, p_{\rm Comp}, \rho e \right\}$:
\begin{equation}
  \frac{\partial \mathbf{Q}}{\partial t} = \left( \begin{array}{c}
   - \mathbf{u} \cdot \nabla \rho - \rho \nabla \cdot \mathbf{u} \\
   - \mathbf{u} \cdot \nabla \mathbf{u} - \frac{1}{\rho} \nabla p_{\rm Comp} \\
   - \mathbf{u} \cdot \nabla p_{\rm Comp} - \rho c^2 \nabla \cdot \mathbf{u} \\
   - \mathbf{u} \cdot \nabla \left( \rho e \right) - \left(\rho e + p_{\rm Comp} \right) \nabla \cdot \mathbf{u}
                     \end{array} \right)
\end{equation}

 \item A piecewise quadratic parabolic profile approximation of $\mathbf{Q}$ is constructed within each cell with a modified version of the PPM algorithm \cite{Almgren:2010a}. These constructions are performed in each coordinate direction separately.
 
 \item Average values of $\mathbf{Q}$ are predicted on edges over the time step using characteristic extrapolation. A characteristic tracing operator with flattening is applied to the integrated quadratic profiles in order to obtain left and right edge states at $k+1/2$
 
 \item  An approximate Riemann problem solver is employed to compute the primitive variables centered in time at $k+1/2$, and in space on the edges of a cell. This state is denoted as the \emph{Godunov state}: 
$\mathbf{Q}^{\rm gdnv} = 
\left\{\rho^{\rm gdnv},   \mathbf{u}^{\rm gdnv}, p_{\rm Comp}^{\rm gdnv}, \left(\rho e \right)^{\rm gdnv} \right\}$.  
The flux vector $\mathbf{F}^{k+1/2}$ can now be constructed and synchronized over all the 
compressible levels involved in the computation. Then, Eq.~(\ref{eqn:compressible_vectorised}) is updated to  $k+1$. 

\end{enumerate}

\subsubsection{Step 3: Computation of compressible elements on the finest compressible level}
\label{subsubsec:step3}

As explained in \S\ref{subsec:overall_presentation}, terms involving the pressure as well as the velocity and its divergence are provided to the low-Mach-number Eqs.~(\ref{eqn:mass_LM}-\ref{eqn:energy_LM}) from the compressible solution so as to retain the acoustic effects. Consequently, several terms on level  $l_{\rm max\_comp}$ have to be computed and interpolated to the low-Mach levels  $\left\{l_{{\rm 
    max\_comp}+1},\cdots,L\right\}$. 

Recall that the evaluation of the velocity field is based on a projection method and requires solution of a variable-coefficient Poisson equation for the pressure. As it will be explained in detail in the following steps, two different velocity fields are involved in the algorithm: a normal velocity located at cell edges and centered in time, and a final state velocity located at cell centers and evaluated at the end of a time-step. Consequently, two different projections are also required right hand sides for these projections will be differently located in both space and time. Similarly, the \emph{acoustic pressure $p_1$} and its gradient will be required at different position in space and time.

The velocity vector and the \emph{acoustic pressure $p_1$} located at time $t^{n+1}$ are obviously taken from the compressible solution computed at the end of the previous step \S\ref{subsubsec:step2}. Note that the 
\emph{acoustic pressure $p_1$} at time $t^{n+1}$ is computed as follows:
\begin{equation}
  p_1^{n+1} = \left(\rho e\right)^{n+1} \left(\gamma -1\right) - p_0
\end{equation}

Following on, the velocity vector and the \emph{acoustic pressure $p_1$} at time $t^{n+1/2}$ are taken from compressible variables at their Godunov state, i.e.  $\mathbf{u}^{\rm gdnv}$ and $p_{\rm Comp}^{\rm gdnv}$, respectively. As the time advancement of the compressible solution may involve several $K$ sub-steps, $\mathbf{u}^{\rm gdnv}$ and $p_{1}^{\rm gdnv}$ are averaged in time as follows:
\begin{align}
\overline{\mathbf{u}^{\rm gdnv}} &= \left(\sum\limits_1^K \mathbf{u}^{\rm gdnv} \right) / K \\
  \overline{p_{1}^{\rm gdnv}} &= \left(\sum\limits_1^K \left( p_{\rm Comp}^{\rm gdnv}-p_0\right)\right) / K  
\end{align}

The gradient terms $\nabla \overline{p_{1}^{\rm gdnv}}, \nabla \overline{\mathbf{u}^{\rm gdnv}}, \nabla 
\mathbf{u}^{n+1}$ are computed on level  $l_{\rm max\_comp}$, and together with $\mathbf{u}^{n+1}$ and $p_1^{n+1}$ are interpolated to the low-Mach levels  $\left\{l_{{\rm max\_comp}+1},\cdots,L\right\}$. Note that except $\nabla \mathbf{u}^{n+1}$ which is nodal, all other terms are located at cell centers.

\subsubsection{Step 4: Computation of material derivative of the acoustic pressure $p_1$}
\label{subsubsec:step4}

The material derivative of the acoustic pressure $p_1$, which appears in 
the RHS of Eq.~(\ref{eqn:energy_LM}), is now computed. This term is computed on all low-Mach levels $\left\{l_{{\rm 
    max\_comp}+1},\cdots,L\right\}$ as follows:
\begin{equation}
  \frac{{\rm D} p_1}{{\rm D} t} = \frac{p_1^{n+1} - p_1^{n}}{\Delta t_{\rm{LM}}} + 
  \frac{\mathbf{u}^{n+1}+\mathbf{u}^{n}}{2} \nabla \overline{p_{1}^{\rm gdnv}} \label{eqn:material_derivative_p1}
\end{equation}
Here, $p_1^{n+1}, \mathbf{u}^{n+1}$ and $\nabla \overline{p_{1}^{\rm gdnv}}$ are already known because they were 
computed 
during the time advancement of the fully compressible Eqs.~(\ref{eqn:mass_comp}-\ref{eqn:energy_comp}) through the 
previous steps described from \S\ref{subsubsec:step2} to \S\ref{subsubsec:step3}. 
Of course, $p_1^{n}$ and $\mathbf{u}^{n}$ are known from the previous time-step iteration.

\subsubsection{Step 5: Time advancement of the low-Mach-number equations: thermodynamic system}
\label{subsubsec:step5}

The low-Mach-number Eqs.~(\ref{eqn:mass_LM}-\ref{eqn:energy_LM}) are now advanced in time on all low-Mach levels, i.e. 
on $\left\{l_{{\rm max\_comp}+1},\cdots,L\right\}$. As explained at the beginning of this section, the set of equations is solved through a fractional step procedure. Consequently, the thermodynamic system composed of Eq.~(\ref{eqn:mass_LM}) and Eq.~(\ref{eqn:energy_LM}) is advanced first. Then the momentum Eq.~(\ref{eqn:momentum_LM}) is advanced with a projection method. The whole procedure is described below. 

The very first step is to compute the normal velocity on the edges of a computational cell and at time $t^{n+1/2}$, which is denoted  $\mathbf{u}^{\rm MAC}$ for convenience. Here the superscript MAC refers to a MAC-type staggered grid \cite{Harlow:1965} discretization of the equations. A provisional value of the normal velocity on edges, denoted $\mathbf{u}^{*,\rm MAC}$, is estimated from $\mathbf{u}^{n}$ with the PPM algorithm. Note that during this procedure, the cell-centered gradients of the pressure, which appear in the RHS of the momentum Eq.~(\ref{eqn:momentum_LM}), are included as an explicit source term contribution for the $1$D characteristic tracing (see ~\cite{Miller:2002}):
\begin{equation}
  S^n =   \frac{1}{\rho^n}\left(\nabla \overline{p_{1}^{\rm 
gdnv}} + \nabla   p_{2}^{n-1/2}\right) 
\end{equation}
Recall here that $\mathbf{u}^{*,\rm MAC}$ is only a provisional value of the normal velocity on edges and a projection operator is applied to ensure that the divergence constraint constructed with the interpolation of $\nabla \overline{\mathbf{u}^{\rm gdnv}}$ is discretely satisfied.

In the numerical resolution of low-Mach-number systems, several different strategies have been developed to ensure the correctness of the solution (see the paper of \citet{Knikker:2011} for a review). Here, the so-called \emph{volume discrepancy} approach is employed. Mass and energy are advanced in a conservative form, however the constraint on the velocity field fails to ensure that the equation of state is satisfied. Thus, the constraint is locally modified by an additional term to maintain a thermodynamic consistency so as to control the drift in pressure from the purely compressible solution. The key observation in volume discrepancy approaches is that local corrections can be added to the constraint in order to specify how the local thermodynamic pressure is allowed to change over a time step to account for the numerical drift. After numerical integration over a time step, for a given cell if the thermodynamic pressure is too low, the net flux into the cell needs to be increased; if it is too high, the net flux needs to be decreased. This is a fundamental concept of volume discrepancy approaches, and a rigorous analysis derived in the context of reactive flows with complex chemistry is given in an upcoming work \cite{Nonaka:2017}.

An iterative procedure is now performed to advance Eq.~(\ref{eqn:mass_LM}) and Eq.~(\ref{eqn:energy_LM}) so as to converge towards a value of $\mathbf{u}^{\rm MAC}$ that ensures the conservation of the equation of state. The provisional velocity $\mathbf{u}^{*,\rm MAC}$ is corrected via a projection method that includes
solution of a variable-coefficient Poisson equation. 
The new value of the velocity is then used to define the convective terms in 
Eq.~(\ref{eqn:mass_LM}) and Eq.~(\ref{eqn:energy_LM}) and to advance 
$\rho$ and $\left(\rho h \right)$. At each iteration, the correction, $\Delta S$, is added to the RHS 
of the Poisson equation so as to control the drift of the low-Mach-number solution from the equation of state given by the fully compressible solution.

Starting from iteration $m=1$, 
\begin{align}
  \nabla \left(\frac{1}{\rho^n} \nabla \phi_m \right) &= \nabla\mathbf{u}^{*,\rm MAC}  - \left( \nabla \overline{\mathbf{u}^{\rm gdnv}} + \Delta S_{m-1}  \right)  \label{eqn:drift_loop_Poisson} \\
\mathbf{u}^{\rm MAC}_{m} &= \mathbf{u}^{*,\rm MAC} - \frac{1}{\rho^n} \nabla \phi_m \label{eqn:advance_umac}\\
\left.\frac{{\rm D} p_1}{{\rm D} t}\right\vert_{m} &= \frac{p_1^{n+1} - p_1^{n}}{\Delta t_{\rm{LM}}} + 
\mathbf{u}^{\rm MAC}_{m} \nabla \overline{p_{1}^{\rm gdnv}} \\
\rho_{m} &= \rho^n - \Delta t_{\rm{LM}} \nabla \left(\mathbf{u}^{\rm MAC}_{m} \rho^{n+1/2} \right) 
\label{eqn:drift_loop_mass_LM} \\
\left(\rho h \right)_{m} &= \left(\rho h \right)^{n} - \Delta t_{\rm{LM}} \nabla \left(\mathbf{u}^{\rm MAC}_{m} 
\left(\rho h\right)^{n + 1/2} \right) + \Delta t_{\rm{LM}} \left.\frac{{\rm D} p_1}{{\rm D} t}\right\vert_{m} 
\label{eqn:drift_loop_energy_LM} 
\end{align}
Here $\rho^{n+1/2}$ and $\left(\rho h\right)^{n+1/2}$ are the edge states predicted with the PPM algorithm from $\rho^n$ and $\left(\rho h\right)^n$, respectively. Note that similarly to the prediction of the velocity $\mathbf{u}^{*,\rm MAC}$ on edges, the cell-centered term ${\rm D} p_1/{\rm D} t$ that appear in the RHS of Eq.~(\ref{eqn:drift_loop_energy_LM}) is taken into account during the computation of  $\left(\rho h\right)^{n+1/2}$ as an explicit source term contribution. Note also that for $m=1$, $\Delta S_{m-1}=0$.

At the end of each iteration, after evaluation of Eq.~(\ref{eqn:drift_loop_energy_LM}), the 
drift in pressure is computed as follows:
\begin{align}
 \delta p_m &= \left(\rho h \right)_{m} \frac{\gamma -1}{\gamma} - \left(p_1^{n+1} + p_0 \right) \\
 \Delta S_{m,i} &= \frac{\delta p_m}{\left(p_1^{n+1}+p_0 \right)\Delta t_{\rm{LM}}} \\
 \Delta S_m &= \Delta S_{m,i} - \frac{1}{V}\int_V \Delta S_{m,i}~{\rm dV} \label{eqn:Delta_S}\\
 \epsilon_m &= \frac{\max\left(|\delta p_m |  \right)}{||p_1^{n+1} + p_0 ||} \label{eqn:epsilon_p}
\end{align}
Here, $|\cdot|$ and $|| \cdot ||$ are the absolute value and the infinity norm, respectively. Note that $\Delta 
S_{m,i}$ 
denotes the point-wise computation of $\Delta S_m$ for each cell $i$. The equation of state is considered satisfied at convergence for 
$\epsilon_m < \epsilon_p$, where $\epsilon_p$ is specified by the user. At convergence, 
$\rho^{n+1} = \rho_{m}$, $\left(\rho h \right)^{n+1}=\left(\rho h \right)_{m}$ and $\mathbf{u}^{\rm MAC} = 
\mathbf{u}^{\rm MAC}_m$.

During this whole procedure, once $\mathbf{u}^{\rm MAC}_{m}$, $\left(\rho h \right)_{m}$, $\left(\rho h \right)_{m}$ and $\Delta S_m$ are evaluated with Eqs.~(\ref{eqn:advance_umac}), (\ref{eqn:drift_loop_mass_LM}), (\ref{eqn:drift_loop_energy_LM}) and (\ref{eqn:Delta_S}), respectively, the variables are synchronized over the levels so as to take into account the contribution of finest levels to the coarser low-Mach-number level $l_{{\rm max\_comp}+1}$.

\subsubsection{Step 6: Time advancement of the low-Mach-number equations: momentum equation}
\label{subsubsec:step6}

The momentum Eq.~(\ref{eqn:momentum_LM}) is now advanced in time with a fractional step, projection 
method. First, a provisional velocity field is computed as follows: 
\begin{align}
  \mathbf{u}^{*,n+1} = \mathbf{u}^n - \Delta t_{\rm{LM}} \left(\overline{\mathbf{u}^{\rm MAC}} \cdot \nabla \mathbf{u}^{n+1/2} \right) - \Delta t_{\rm{LM}} \left(\frac{1}{\rho^{n+1/2}}\nabla \overline{p_{1}^{\rm gdnv}} + 
  \frac{1}{\rho^{n+1/2}} \nabla p_2^{n-1/2} \right)
\end{align}
with $\rho^{n+1/2}=\left(\rho^{n+1}+\rho^{n} \right)/2$. Recall that $\mathbf{u}^{\rm MAC}$ lives on the edges of a computational cell, $\overline{\mathbf{u}^{\rm MAC}}$  represents the spatial average to cell centers. 
Again, $\mathbf{u}^{n+1/2}$ is the prediction of the time and space centered values of the velocity $\mathbf{u}^{n}$ via the PPM algorithm, and the terms $\left(\frac{1}{\rho^{n}}\nabla \overline{p_{1}^{\rm gdnv}} + 
  \frac{1}{\rho^{n}} \nabla p_2^{n-1/2} \right)$ are taken into account during the construction of $\mathbf{u}^{n+1/2}$ as an explicit source term contribution.

The following variable-coefficient Poisson equation for the pressure is solved to enforce the divergence 
constraint on the velocity field:
\begin{equation}
  \nabla \cdot \left( \frac{1}{\rho^{n+1/2}} \nabla \phi\right) = \nabla \cdot \left(\mathbf{u}^{*,n+1} + \frac{ \Delta 
    t_{\rm{LM}}}{\rho^{n+1/2}}\nabla p_2^{n-1/2} \right) - \left.\left(\nabla \cdot \mathbf{u}^{n+1}\right)\right\vert_{\rm Comp} \label{eqn:poisson_eq}
\end{equation}
Note that a subscript $\rm Comp$ has been added here to $\nabla \mathbf{u}^{n+1}$ in order to recall that it has been computed from 
the solution of the fully compressible equations and has been interpolated from the compressible level 
$\left\{l_{\rm max\_comp}\right\}$ to the low-Mach levels $\left\{l_{{\rm max\_comp}+1},\cdots,L\right\}$.

Finally, the provisional velocity field $\mathbf{u}^{*,n+1}$ is corrected as follows:
\begin{equation}
  \mathbf{u}^{n+1} = \mathbf{u}^{*,n+1} - \frac{1}{\rho^{n+1/2}} \nabla \phi
  \label{eqn:advance_u_np1}
\end{equation}
and the hydrodynamic pressure is also updated:
\begin{align}
  p_2^{n+1/2} &= \frac{1}{\Delta t_{\rm{LM}}} \phi \label{eqn:advance_p2} \\
  \nabla p_2^{n+1/2} &= \frac{1}{\Delta t_{\rm{LM}}} \nabla\phi   \label{eqn:advance_grad_p2}
\end{align}

Similarly to \S\ref{subsubsec:step5}, once $\mathbf{u}^{n+1}$, $p_2^{n+1/2}$ and $\nabla p_2^{n+1/2}$ are evaluated with Eqs.~(\ref{eqn:advance_u_np1}), (\ref{eqn:advance_p2}) and (\ref{eqn:advance_grad_p2}), respectively, the variables are synchronized over the levels so as to take into account the contribution of finest levels to the coarser low-Mach-number level $l_{{\rm max\_comp}+1}$.

\subsubsection{Step 7: Synchronization between the low-Mach-number system and the fully compressible system.}
\label{subsubsec:step7}

The variables $\rho^{n+1}$, $\left(\rho h \right)^{n+1}$ and $\mathbf{u}^{n+1}$ computed on the low-Mach level $l_{{\rm max\_comp}+1}$ are restricted back on the set of compressible levels $\left\{1,\cdots,l_{\rm max\_comp}\right\}$. This operation sets coarse cell-centered values equal to the average of the fine cells covering it. The conservative state variables in Eq.~(\ref{eqn:vector_conservative_vars}) are then updated to take into account the low-Mach-number contribution as follows:

\begin{equation}
  \mathbf{U}^{n+1}=\left\{ \begin{array}{c}
                      \rho^{n+1} \\
		     \rho^{n+1} \mathbf{u}^{n+1} \\
		      \left(\rho e\right)^{n+1} +  \frac{1}{2}  \rho^{n+1} \mathbf{u}^{n+1} \cdot \mathbf{u}^{n+1}
                     \end{array} \right\}
\end{equation}
with $\left(\rho e\right)^{n+1} = \left(\rho h \right)^{n+1} / \gamma$. Of course, this update of the conservative variables is only performed in regions where compressible levels lie beneath low-Mach-number levels. 

Finally, the computation through the time-step is finished and the next iteration can begin at \S\ref{subsubsec:step1}.

\section{Results}
\label{sec:results}

The performance of the new hybrid compressible/low-Mach method proposed in the present paper is now assessed with several test cases. The first test case consists of the propagation of uni-dimensional acoustic waves. The goal of this canonical simulation is to assess the spatial and temporal rates of convergence of the hybrid method. 
The second test case consists of the simultaneous propagation of mixed acoustic, entropic and vorticity modes in a 2D square domain. Finally, a more practical problem similar to the ones encountered in the industry is investigated by simulating the propagation of aeroacoustic waves generated by the formation of a Kelvin-Helmholtz instability in mixing layers. A 
feature of this problem is that a very fine discretization of the mixing layer interface is required to accurately capture 
the vortex formation.
It will be demonstrated that in the context of an AMR framework, the hybrid method proposed in the present paper leads to larger time-steps by solving the low-Mach-number equations instead 
of the purely compressible equations in the finest levels of discretization.

\subsection{1D acoustic wave propagation}
\label{subsec:results_1d_acoustic}

The first test case consists of the simulation of uni-dimensional acoustic wave propagation in a fluid at rest. 
The computational domain is a rectangle of length $L_x=1$~m and height $L_y=0.125$~m, so that the velocity vector contains only two components $u_x$ and $u_y$, and is periodic in both directions. 
The initial conditions are given as
%
\begin{align}
\rho^{\rm init}\left(x\right) &= \rho_{\rm ref} + A \exp \left(- \left(\frac{x- L_x/2}{\sigma} \right)^2\right) \label{eqn:1d_entropy_pulse_rho} \\
u^{\rm init}_x\left(x\right) &= 0, \hspace{0.5cm} u^{\rm init}_y\left(x\right) = 0 \\
p_0^{\rm init}\left(x\right) &= p_{\rm ref}, \hspace{0.5cm} p_1^{\rm init}\left(x\right) = \rho^{\rm init}\left(x\right) c_0^2
\end{align}
with $A=0.1$ and $\sigma=0.1,$ a set of parameters designed to control the amplification and the width of the acoustic pulse, 
respectively,  while $\rho_{\rm ref}=1.4$~kg/m$^3$, $p_{\rm ref}=10000$~Pa and $c_0$ the sound speed defined as $c_0 = \sqrt{\gamma p_{\rm ref}/\rho_{\rm ref}}=100$~m/s. The heat capacity ratio is set to $\gamma = 1.4$, while the tolerance parameter $\epsilon_p$ in Eq.~(\ref{eqn:epsilon_p}) is set to $\epsilon_p = 1 \times 10^{-13}$ to ensure that no errors are introduced by the drift in pressure of the low-Mach-number solution within the hybrid algorithm. The simulations are performed over $1 \times 10^{-2}$~s, so that $2$ acoustic waves travels through the computational domain in the left and right direction from the initial pulse, and then merge at the end of the simulation to form the same shape as the initial pulse.

Consider a simulation with 6 levels, and define
$N^l_x$ and $N^l_y$ as the number of cells at level $l$ in the $x$ and $y$ directions, respectively.
The first level $l=1$ is discretized with $N^{l=1}_x=32$ and $N^{l=1}_y=4$ points, 
while the other levels are progressively discretized with a mesh refinement ratio of a factor of $2$. Note here that the whole domain is covered by all the levels. Table~\ref{tab:summary_1D_acoustic_meshgrid} summarizes the configuration.

\begin{table}
\centering
\renewcommand\arraystretch{1.5}
\begin{tabular}{c || c | c | c | c | c | c }
\renewcommand*{\arraystretch}{0.5}
$l$    & $1$  & $2$ & $3$ & $4$ & $5$ & $6$  \\ \hline\hline 
$N_x$  & $32$ & $64$ & $128$ & $256$ & $512$ & $1024$  \\\hline
$N_y$  & $4$  & $8$  & $16$ & $32$ & $64$ & $128$   
\end{tabular}
\caption{Summary of the configuration for simulations performed on the 1D acoustic waves propagation test case.}
\label{tab:summary_1D_acoustic_meshgrid}
\end{table} 

For all the simulations, the fully compressible Eqs.~(\ref{eqn:mass_comp}-\ref{eqn:energy_comp}) are solved only on one selected level $l_{\rm Comp}=l$. The procedures to perform convergence tests are as follows:
\begin{itemize}
\item for the spatial accuracy, simulations are performed by first selecting, between $l=1$ to $l=5,$ 
the level $l_{\rm Comp}$ where the fully compressible Eqs.~(\ref{eqn:mass_comp}-\ref{eqn:energy_comp}) are solved, and then by selecting a successive addition of low-Mach-number levels of mesh refinement, the finest level chosen being designed by $L$. In total, $15$ simulations are performed, and the choices of $l_{\rm Comp}$ and $L$ for each simulation are summarized in Table~\ref{tab:summary_1D_acoustic}. Furthermore, the low-Mach-number time-step $\Delta t_{\rm{LM}}$ is kept at $9.0 \times 10^{-5}$~s, which corresponds to the minimum time-step for a CFL condition $\sigma^{{\rm CFL}} =0.5$ and for the finest level of refinement $L=6$. Consequently, for all simulations $\Delta t_{\rm{Hyb}}$ is equal to $\Delta t_{\rm{LM}}$ and $K=1$.

\begin{table}
\centering
\renewcommand\arraystretch{1.5}
\begin{tabular}{c || c | c | c | c | c  | c }
\renewcommand*{\arraystretch}{0.5}
\diagbox{$l_{\rm Comp}$}{$L$} & $1$ & $2$ & $3$ & $4$ & $5$ & $6$  \\ \hline\hline 
$1$ &  & $\times$ & $\times$ & $\times$ &  $\times$ &  $\times$ \\\hline
$2$ &  &   & $\times$ & $\times$ &  $\times$ &  $\times$ \\\hline
$3$ &  &  &  & $\times$ &  $\times$ &  $\times$ \\\hline
$4$ &  &  &  &  &  $\times$ &  $\times$ \\\hline
$5$ &  &  &  &  &  &  $\times$ \\
\end{tabular}
\caption{Summary of the choices of $l_{\rm Comp}$ and $L$ for all simulations performed during spatial convergence test of the hybrid method with the propagation of a uni-dimensional acoustic wave.}
\label{tab:summary_1D_acoustic}
\end{table}

\item for the temporal accuracy, the fully compressible Eqs.~(\ref{eqn:mass_comp})-(\ref{eqn:energy_comp}) are solved 
on $l_{\rm Comp}=5$, while the low-Mach-number Eqs.~(\ref{eqn:mass_LM})-(\ref{eqn:energy_LM}) are solved on the last and finest level of mesh refinement $L=6$ so as to minimize spatial discretization errors. Simulations are performed with successive time-steps of $\Delta t_{\rm{LM}}= 3.125, 6.25, 12.5, 25, 50, 100, 200, 400  \times 10^{-6}$~s. Note that for all simulations, the compressible time-step $\Delta t_{\rm Comp}$ is not imposed but computed with Eqs.~(\ref{eqn:delta_t_hyb}) and (\ref{eqn:ceiling_time_steps}).
\end{itemize}

Convergence tests are evaluated with the $\mathcal{L}^2$-norm of the difference on the density between the computed and the initial solution defined by Eq.~(\ref{eqn:1d_entropy_pulse_rho}), which is expressed as follows:
\begin{equation}
\varepsilon_{\rho} = \mathcal{L}^2_{\rho} \left(S_{sol} - S_{init} \right) = \sqrt{\frac{\left(\rho_{sol} - \rho_{init} \right)^2}{N_x^L}}
\end{equation} 
where subscripts $sol$ and $init$ identify the computed and initial solutions $S$. Note that the finest level $L$ of mesh refinement is chosen to compare with the initial solution.

\begin{figure}[!ht]
\centering
\includegraphics[width=0.65\textwidth]{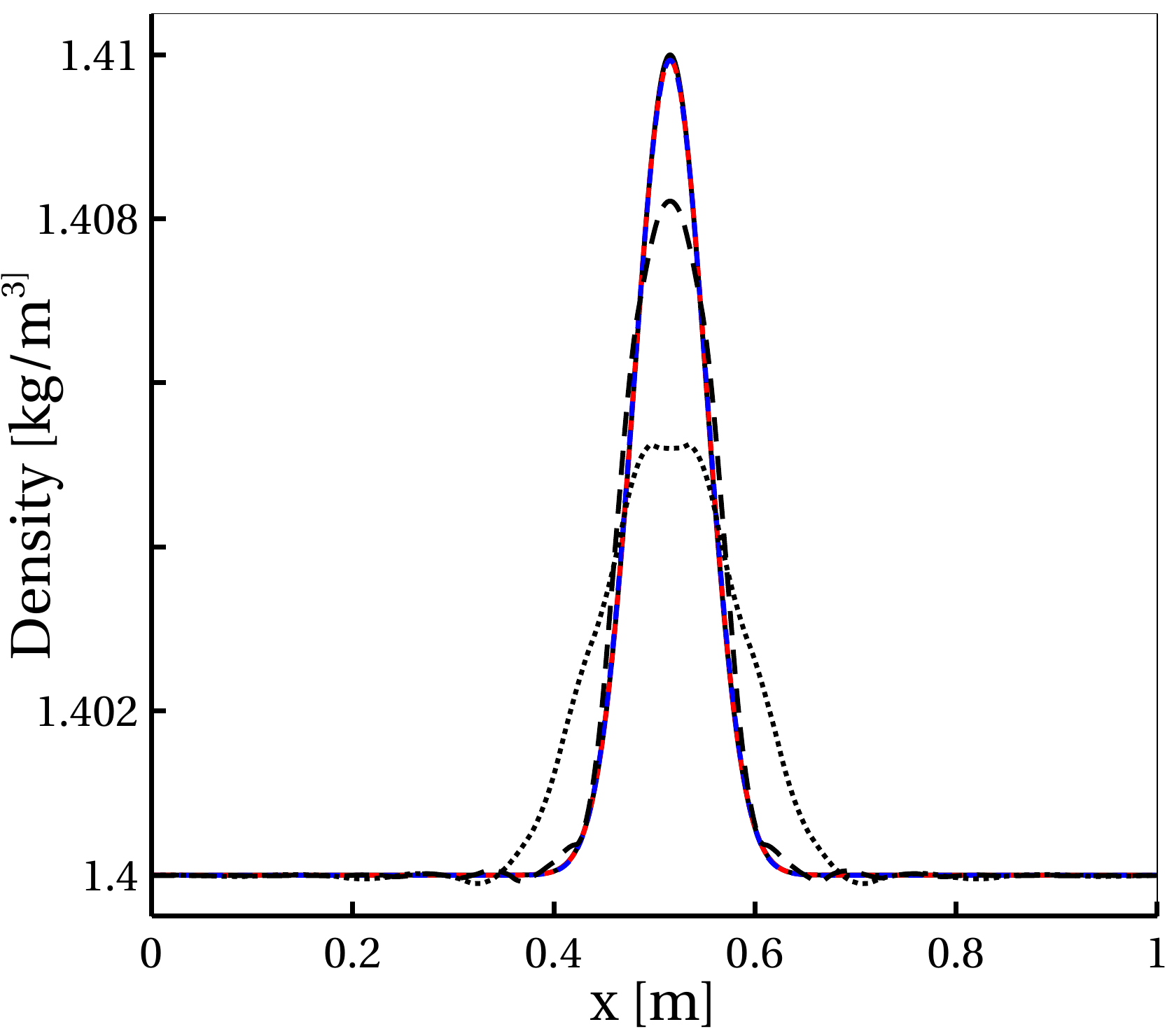}
 \caption{Density profile along $x-$axis. Solid black line: initial acoustic pulse at $0$~s. Computed solutions with $N_x^{l_{\rm Comp}}=32$ (black dotted line), $N_x^{l_{\rm Comp}}=64$ (black dashed line), $N_x^{l_{\rm Comp}}=128$ (blue dotted line) and $N_x^{l_{\rm Comp}}=256$ (red dashed line) at $1 \times 10^{-2}$~s after the merge of the two traveling acoustic waves.}
 \label{fig:results_1D_acoustic_wave_visu_profile}
\end{figure}

\begin{figure}[!ht]
\centering
\includegraphics[width=0.65\textwidth]{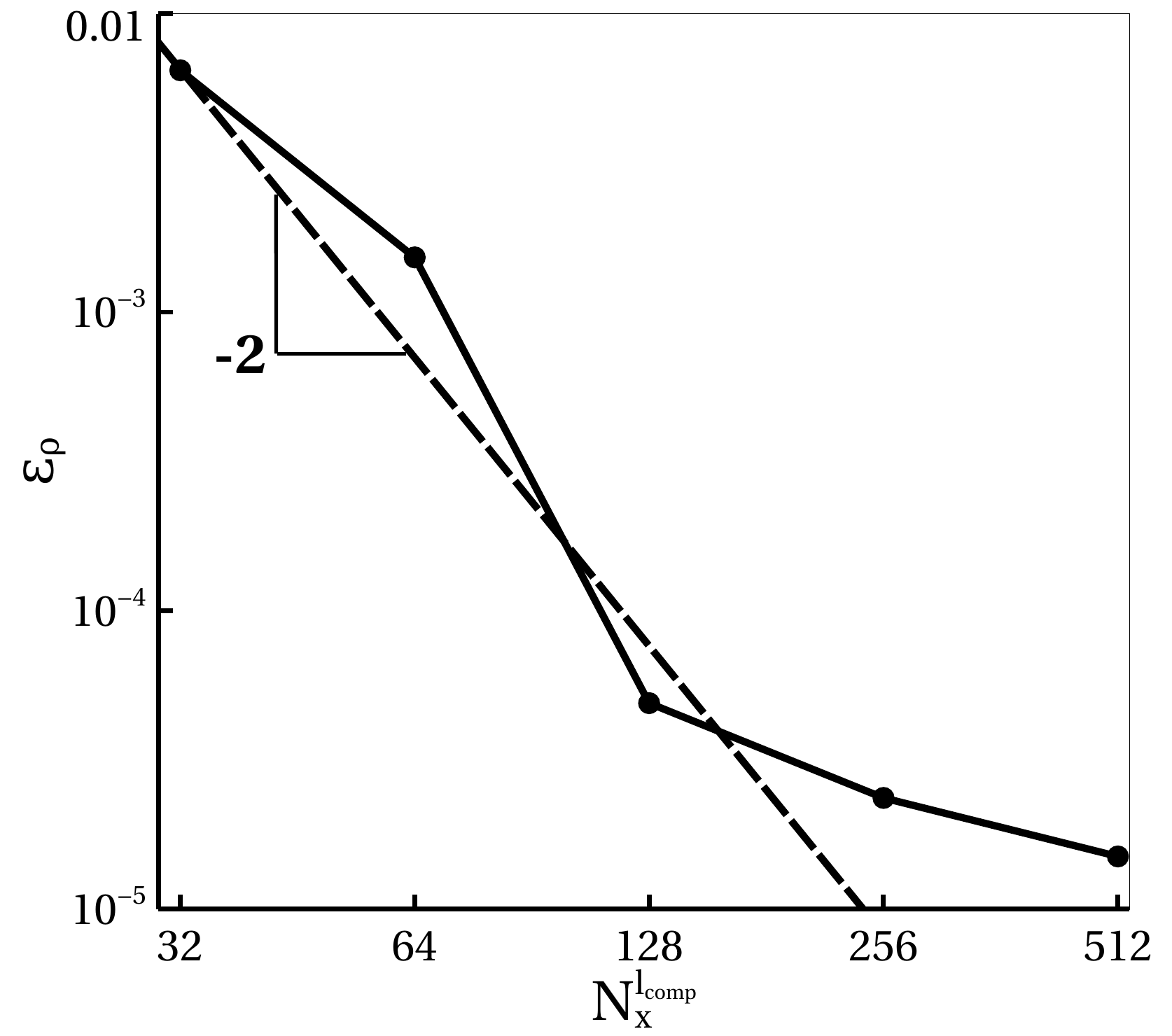}
 \caption{$\mathcal{L}^2$-norm of the discretization error $\varepsilon_{\rho}$ computed for the density, with $L=6$ ($N_x^L=1024$) and $l_{\rm Comp}$ set at different levels $l=1$ to $5$ ($N_x^{l_{\rm Comp}}=32$ to $N_x^{l_{\rm Comp}}=512$). The dashed black line represent a second order slope.}
 \label{fig:results_1D_acoustic_wave_convergence_spatial}
\end{figure}

Figures~\ref{fig:results_1D_acoustic_wave_visu_profile} and~\ref{fig:results_1D_acoustic_wave_convergence_spatial} present profiles of the density as well as the discretization error $\varepsilon_{\rho}$, respectively, for $L=6$ ($N_x^L=1024$) and $l_{\rm Comp}$ set at different levels $l=1$ to $5$ ($N_x^{l_{\rm Comp}}=32$ to $N_x^{l_{\rm Comp}}=512$). In Figure~\ref{fig:results_1D_acoustic_wave_visu_profile} it is observed that under-resolution of the mesh leads to significant dissipation and dispersion of the acoustic waves. Note that the solution computed with $N_x^{l_{\rm Comp}}=512$ is virtually similar to the one computed with $N_x^{l_{\rm Comp}}=256$, and thus is not displayed for clarity purpose. The discretization error $\varepsilon_{\rho}$ is reported in Figure~\ref{fig:results_1D_acoustic_wave_convergence_spatial}, and it is observed that $\varepsilon_{\rho}$ follows a global convergence rate of second-order, which was expected because the algorithm employs a second-order Godunov procedure. Moreover, it can be seen that for $N_x^{l_{\rm Comp}} > 128$, the error starts to reach a plateau with a first-order behavior. This can be explained by the fact that from $32 < N_x^{l_{\rm Comp}} < 128$ the error is dominated by the resolution on the compressible grids, hence a second order accuracy resulting from the second order Godunov method is seen.  At higher resolutions the compressible solution is sufficiently accurate that the error measured is a combination of that from the compressible and low-Mach-number grids, which results in the apparent reduction in order because in this study the low-Mach-number resolution does not change.

\begin{figure}[!ht]
\centering
\includegraphics[width=0.65\textwidth]{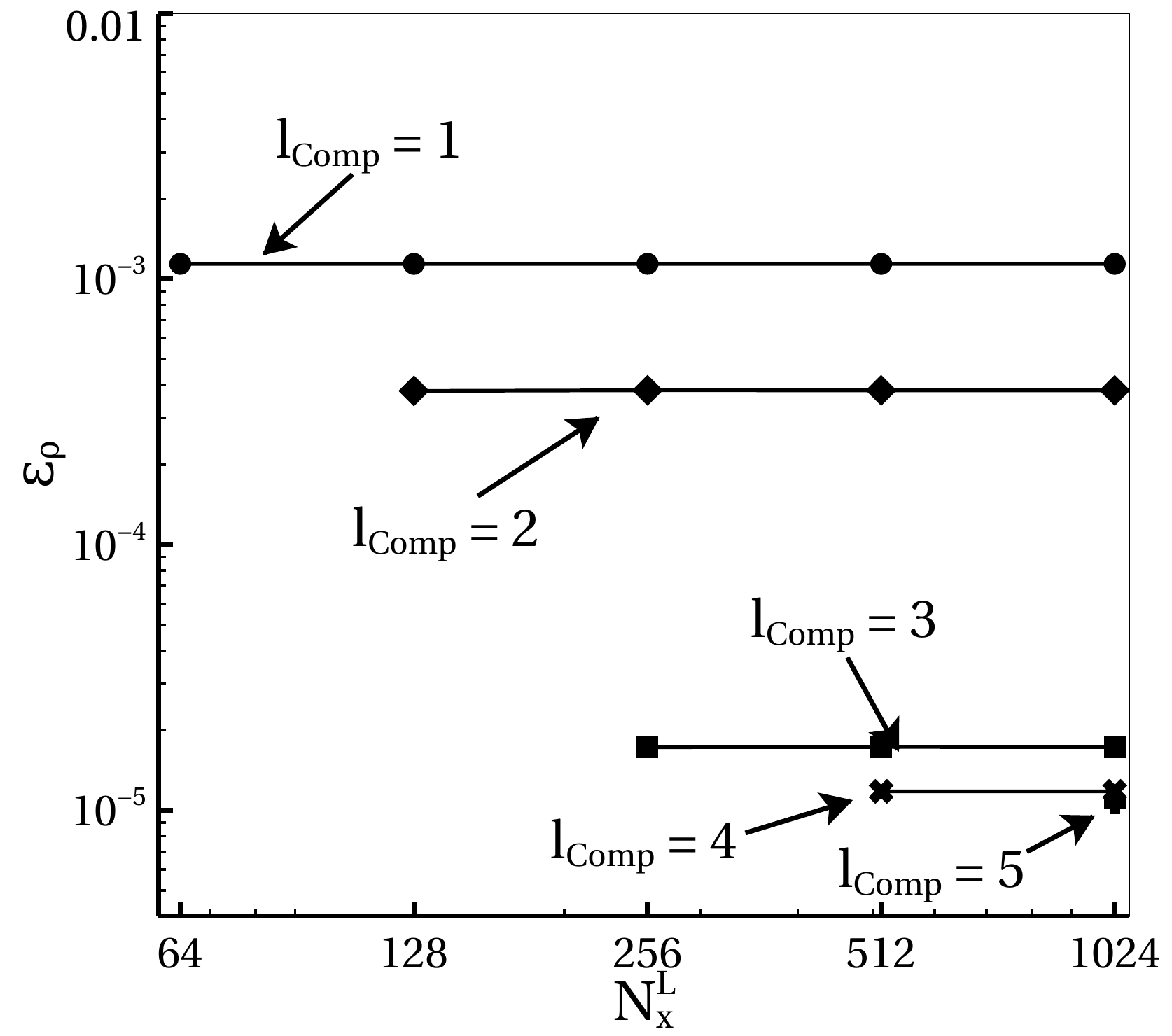}
 \caption{$\mathcal{L}^2$-norm of the discretization error $\varepsilon_{\rho}$ computed for the density and for different maximum level of mesh refinement $L$ where the low-Mach-number equations are solved. Circle, diamond, square, cross and plus symbols represent the fully compressible equations solved on the level $l_{\rm Comp}$ set at $l=1$, $l=2$, $l=3$, $l=4$ and $l=5$, respectively.}
 \label{fig:results_1D_acoustic_wave_convergence_spatial_for_LM}
\end{figure}

The effect of solving the low-Mach-number equations on additional levels of mesh refinement, and for $l_{\rm Comp}$ set at different levels, is shown in Figure~\ref{fig:results_1D_acoustic_wave_convergence_spatial_for_LM}. Circle, diamond, square, cross and plus symbols represent $l_{\rm Comp}$ set at $l=1$, $l=2$, $l=3$, $l=4$ and $l=5$, respectively. This corresponds to a discretization of $N_x^{l_{\rm Comp}}=32, 64, 128, 256$ and $512$ points, respectively. As reported above, the discretization error $\varepsilon_{\rho}$ is reduced as the compressible equations are solved on the finest level. 
In contrary, solving the low-Mach-number equations on finer levels of mesh refinement has no impact on the solution. This behavior was expected, because as the simulation involves only a purely acoustic phenomenon, it is emphasized that the contribution of the set of low-Mach-number equations should be negligible.

\begin{figure}[!ht]
\centering
\includegraphics[width=0.65\textwidth]{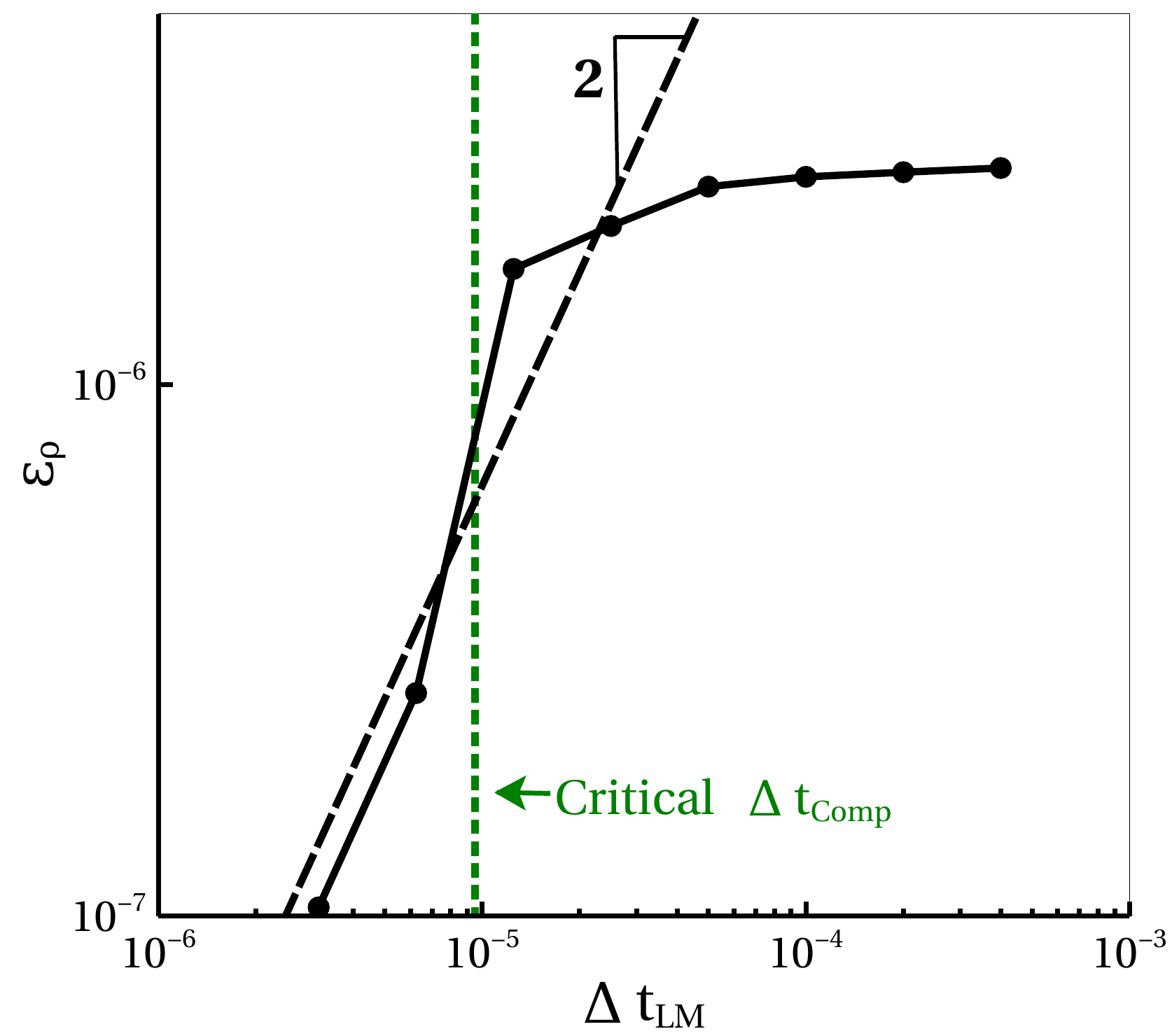}
 \caption{$\mathcal{L}^2$-norm of the discretization error $\varepsilon_{\rho}$ computed for the density for different values of $\Delta t_{\rm LM}$, and with  $l_{\rm Comp}=5$ ($N_x^{l_{\rm Comp}} = 512$) and $L=6$ ($N_x^L=1024$). The dashed black line represent a second order slope.}
 \label{fig:results_1D_acoustic_wave_convergence_time}
\end{figure}

Figure~\ref{fig:results_1D_acoustic_wave_convergence_time} presents the discretization error $\varepsilon_{\rho}$ for different values of $\Delta t_{\rm LM}$. Recall that for these simulations  $l_{\rm Comp}=5$ ($N_x^{l_{\rm Comp}} = 512$) and $L=6$ ($N_x^L=1024$), the corresponding maximum critical compressible time-step for stability and for a CFL condition $\sigma^{{\rm CFL}} =0.5$ is approximately $\Delta t_{\rm Comp}^{\rm crit}=9.5 \times 10^{-6}$~s and is represented in Figure~\ref{fig:results_1D_acoustic_wave_convergence_time} by the dashed vertical green line. It is interesting to notice that when $\Delta t_{\rm LM}$ is larger than the critical time-step, $\Delta t_{\rm Hyb}$ is always set to $\Delta t_{\rm Comp}^{\rm crit}$ and the convergence rate is very low. This makes sense, because as the test case features only purely acoustic phenomena, the set of compressible equations dominate the solution. Consequently, for $\Delta t_{\rm LM}> \Delta t_{\rm Comp}^{\rm crit}$ the compressible equations are always advanced with the same compressible time-step within one low-Mach time-step, and only the number of sub-iterations $K$ will change. In contrary, when $\Delta t_{\rm LM}$ becomes smaller than $\Delta t_{\rm Comp}^{\rm crit}$, $\Delta t_{\rm Hyb}=\Delta t_{\rm LM}$ and a second-order convergence rate in time becomes observable.

The convergence studies performed highlight that care must be taken with the hybrid method. It demonstrates that solving the low-Mach-number equations on additional level of mesh refinement is useless on purely acoustic phenomena, and that the proper resolution of the acoustics has a limiting effect on the accuracy of the solution and the performance of the method. In order to investigate more closely this numerical behavior, a more complex test case involving different mixed modes of fluctuations is now computed, with a solution being a combination of purely acoustics propagation, purely entropic and vorticity convection.

\subsection{2D mixed waves propagation}
\label{subsec:results_2d_tamwebb}

The present test case consists of the propagation and convection of mixed acoustic, entropic and vorticity modes in a 2D square domain \cite{Tam:1993}. A mean flow is imposed throughout the domain, and an acoustic pulse is placed in the center of the domain, while entropy and vorticity pulses are initialized downstream. These latter pulses are simply convected by the mean flow, while the acoustic pulse generates a circular acoustic wave which radiates throughout the domain in all directions.  Furthermore, non-reflecting outflow boundary conditions are imposed in all directions of the domain using the Ghost Cells Navier Stokes Characteristic Boundary Conditions (GC-NSCBC) method \cite{Motheau:2017aa}.

The initial conditions are imposed as follows:
\begin{align}
\rho^{\rm init}\left(x,y\right) &= \rho_{\rm ref} + \eta_a e^{- \alpha_a \left(\left(x-x_a \right)^2 + \left(y-y_a \right)^2\right)} +  \eta_e e^{- \alpha_e \left(\left(x-x_e \right)^2 + \left(y-y_e \right)^2\right)} \\
u^{\rm init}\left(x,y\right) &= M c_{\rm ref} + \left(y-y_v\right) \eta_v  e^{- \alpha_v \left(\left(x-x_v \right)^2 + \left(y-y_v \right)^2\right)} \\ 
v^{\rm init}\left(x,y\right) &= -\left(x-x_v\right) \eta_v  e^{- \alpha_v \left(\left(x-x_v \right)^2 + \left(y-y_v \right)^2\right)} \\
p_0^{\rm init}\left(x,y\right) &= \frac{c_{\rm ref}^2 \rho_{\rm ref}}{\gamma}, \hspace{0.5cm} p_1^{\rm init}\left(x,y\right) = c_{\rm ref}^2 \eta_a e^{- \alpha_a \left(\left(x-x_a \right)^2 + \left(y-y_a \right)^2\right)}
\end{align}
%
Here the sound speed $c_{\rm ref}=200$~m/s and the Mach number $M=0.2$, with $\gamma=1.1$ and density $\rho_{\rm ref} = 1$~kg/m$^3$. 
The domain is a square with sides of length $L_x=L_y=256$~m. 
In the above expressions, $\alpha_x$ is related to the semi-length of the Gaussian $b_x$ by the relation 
$\alpha_x = \ln 2/b_x^2$. Finally, the strengths of the pulses are controlled by the following set of parameters:
\begin{align}
b_a = 15, \hspace{.5cm} \eta_a = 0.001, & \hspace{.5cm} x_a = L_{x}/2, \hspace{.62cm} y_a = L_{x}/2 \\
b_e = 5, \hspace{.5cm} \eta_e = 0.0001, & \hspace{.5cm} x_e = 3L_{x}/4, \hspace{.5cm} y_e = L_{x}/2 \\
b_v = 5, \hspace{.5cm} \eta_v = 0.0004, & \hspace{.5cm} x_v = 3L_{x}/4,\hspace{.5cm} y_v = L_{x}/2 
\end{align}

The test case is computed with $3$ different approaches:  
\begin{itemize}
	\item the new hybrid method developed in the present paper,
	\item by solving only the purely low-Mach-number equations (see Sec.~\ref{subsec:eqs_on_LM_layers}),
	\item by solving only the purely compressible equations (see Sec.~\ref{subsec:eqs_on_compressible_layers}).
\end{itemize}

\begin{figure}[ht]
 \begin{subfigmatrix}{4}
  \subfigure[Time $0.1$~s]{\includegraphics{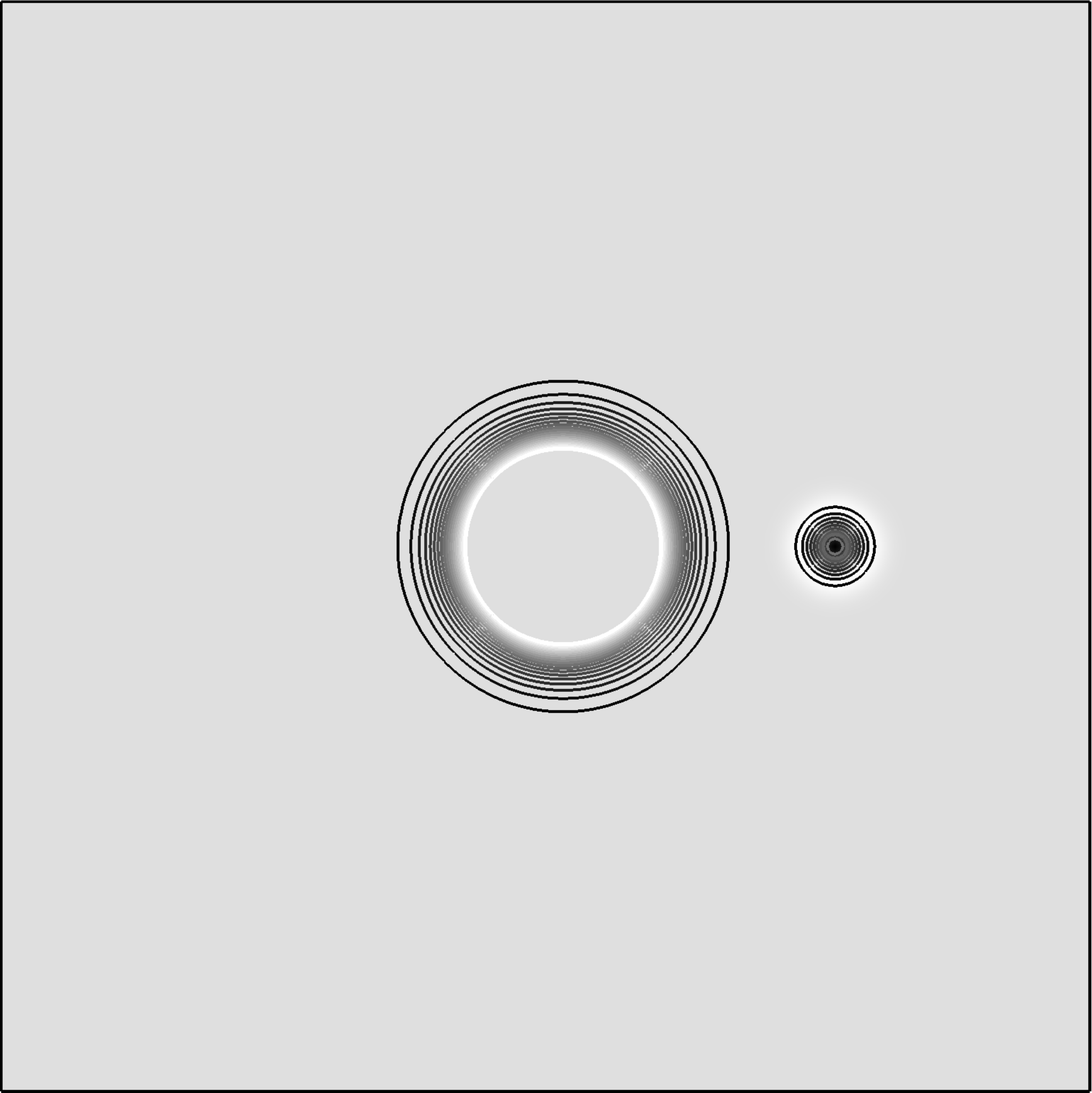}}
  \subfigure[Time $0.2$~s]{\includegraphics{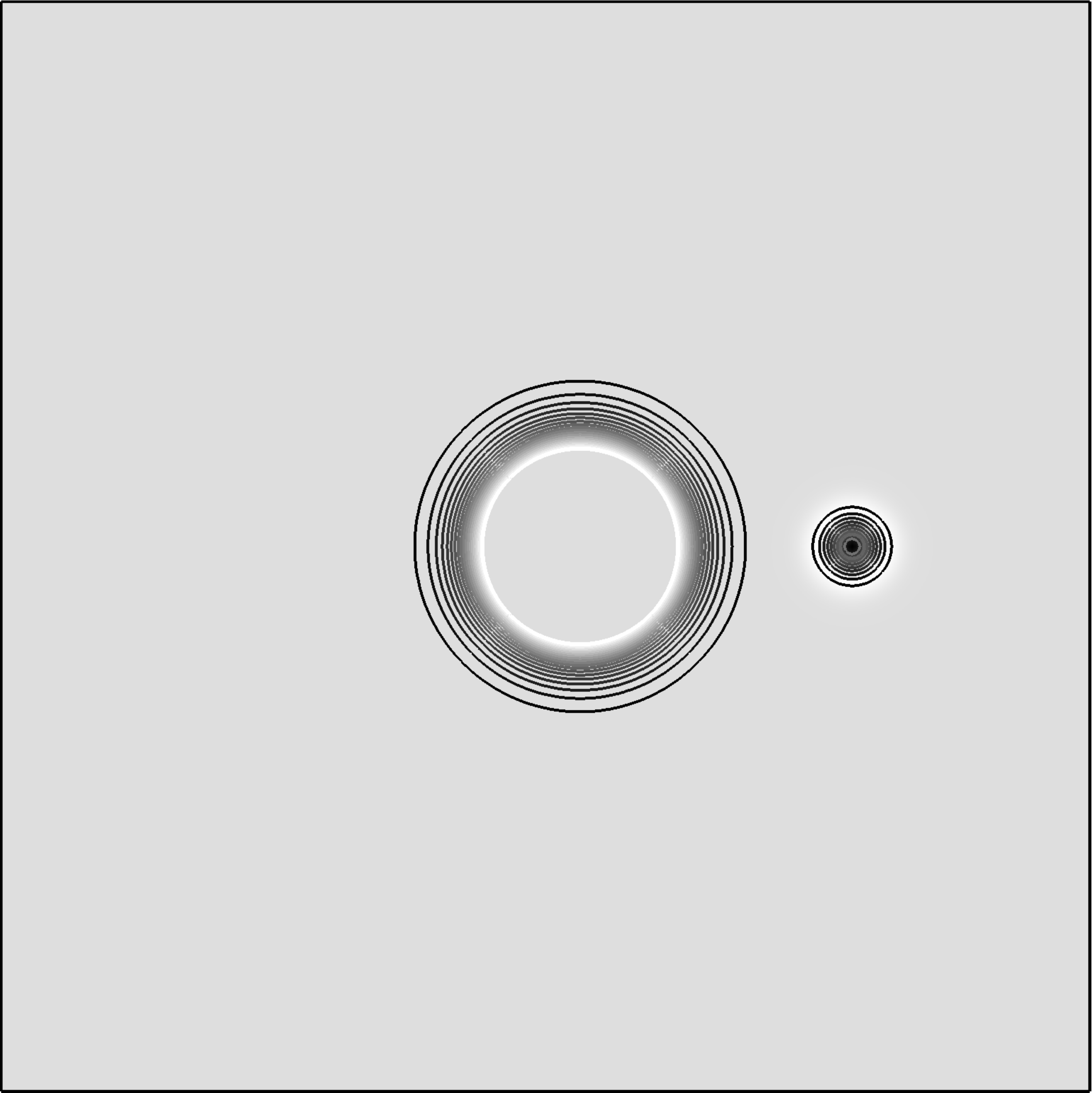}}
  \subfigure[Time $0.3$~s]{\includegraphics{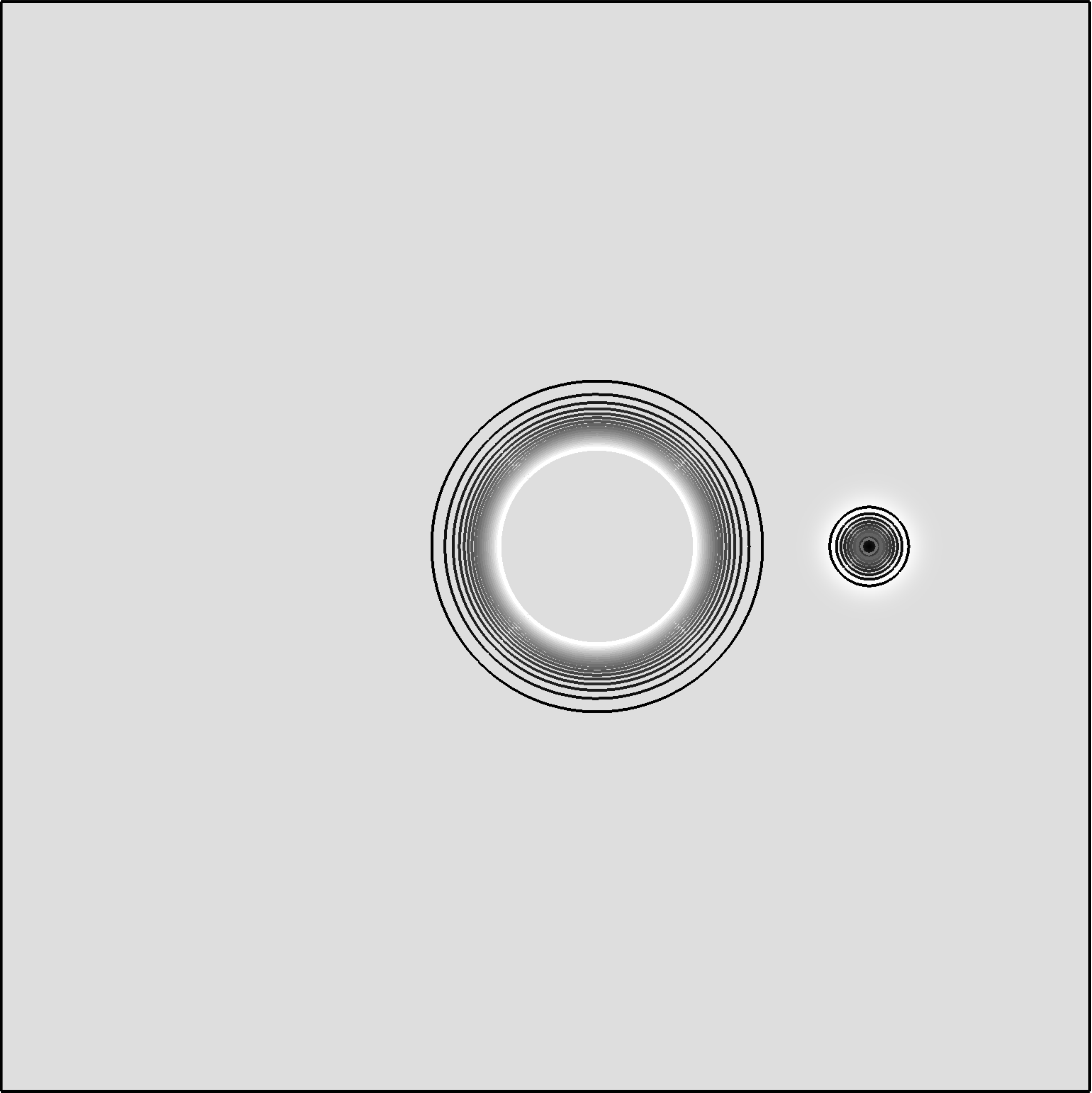}}
  \subfigure[Time $0.4$~s]{\includegraphics{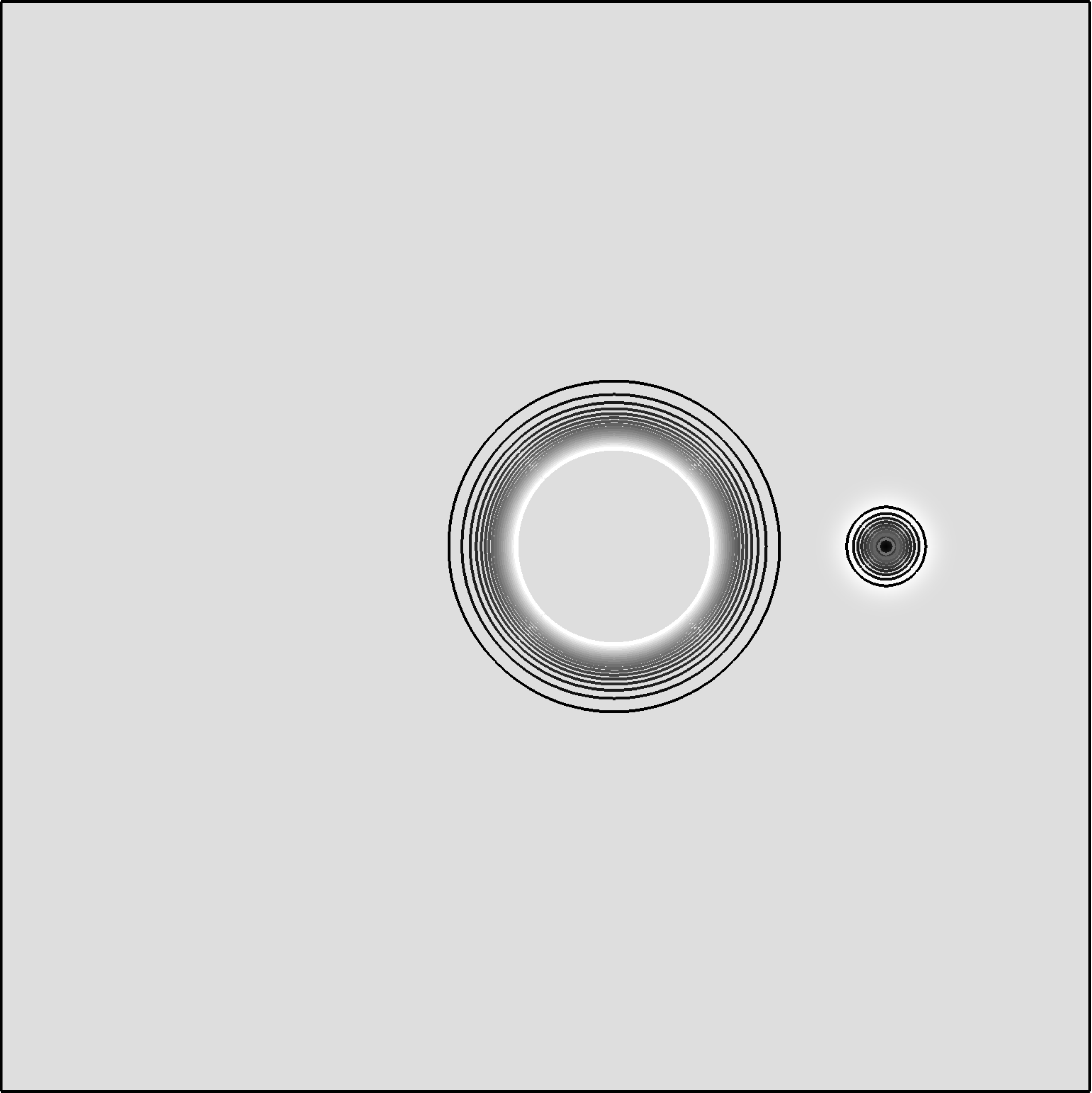}}
  \subfigure[Time $0.1$~s]{\includegraphics{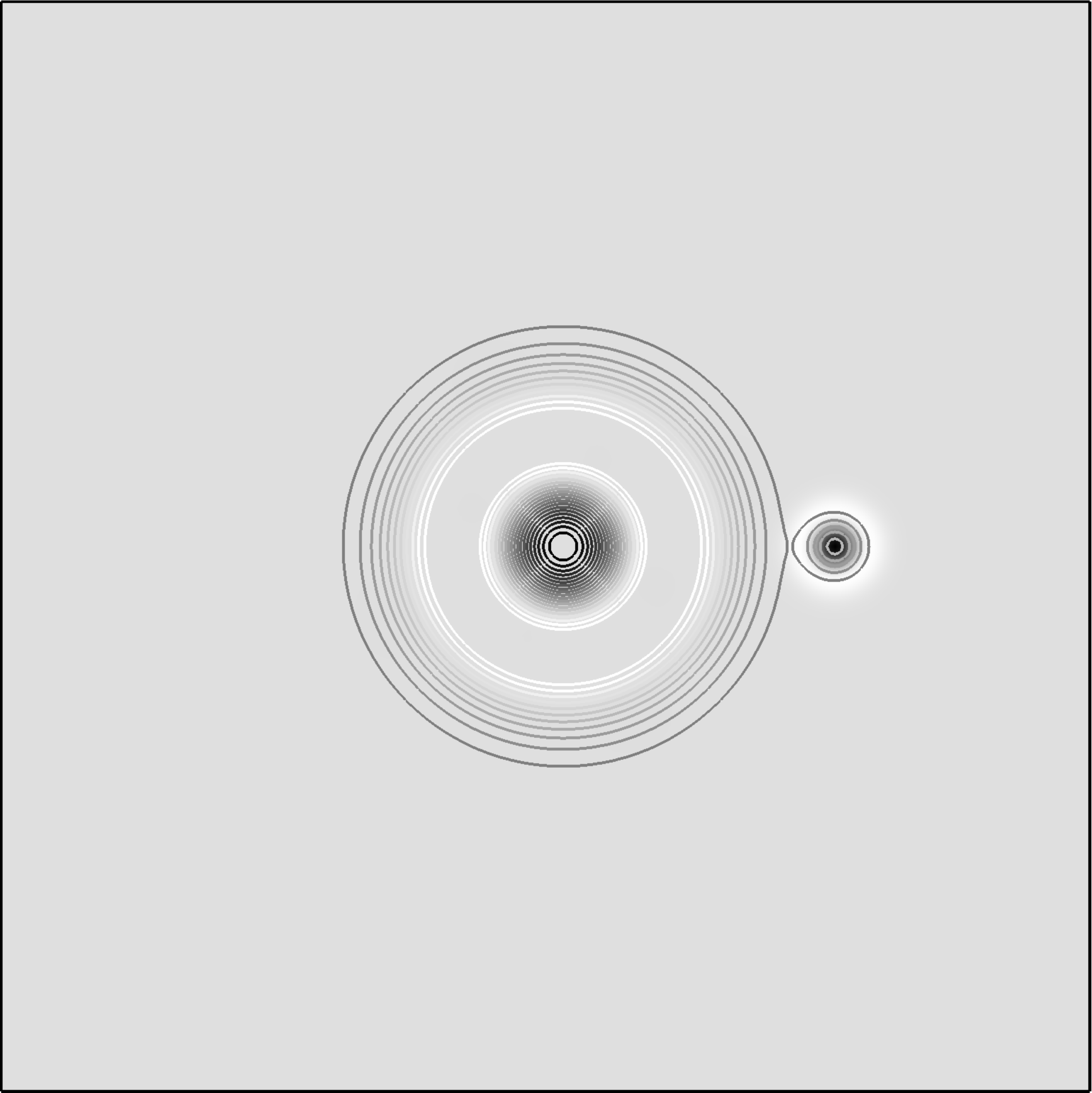}}
  \subfigure[Time $0.2$~s]{\includegraphics{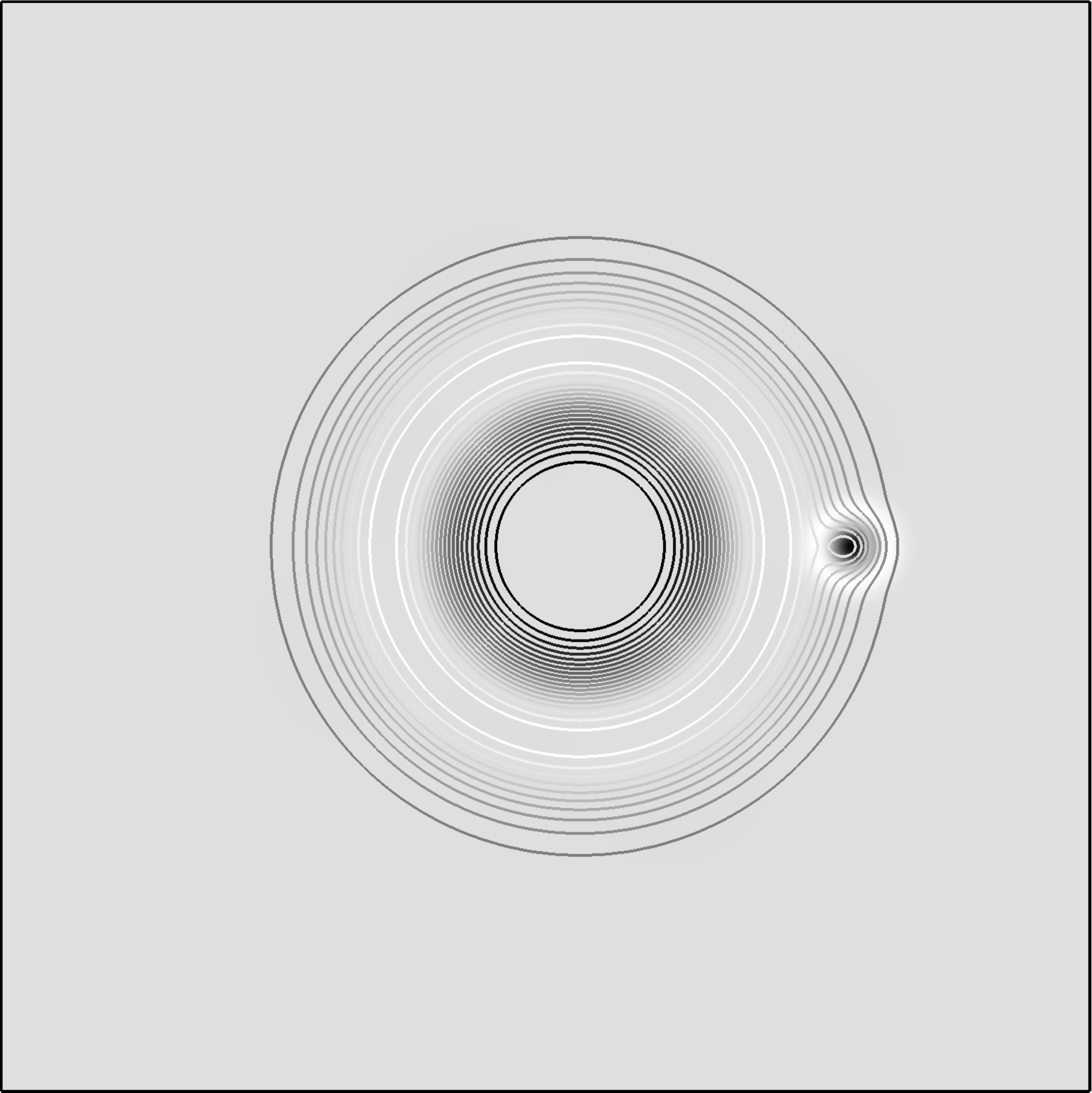}}
  \subfigure[Time $0.3$~s]{\includegraphics{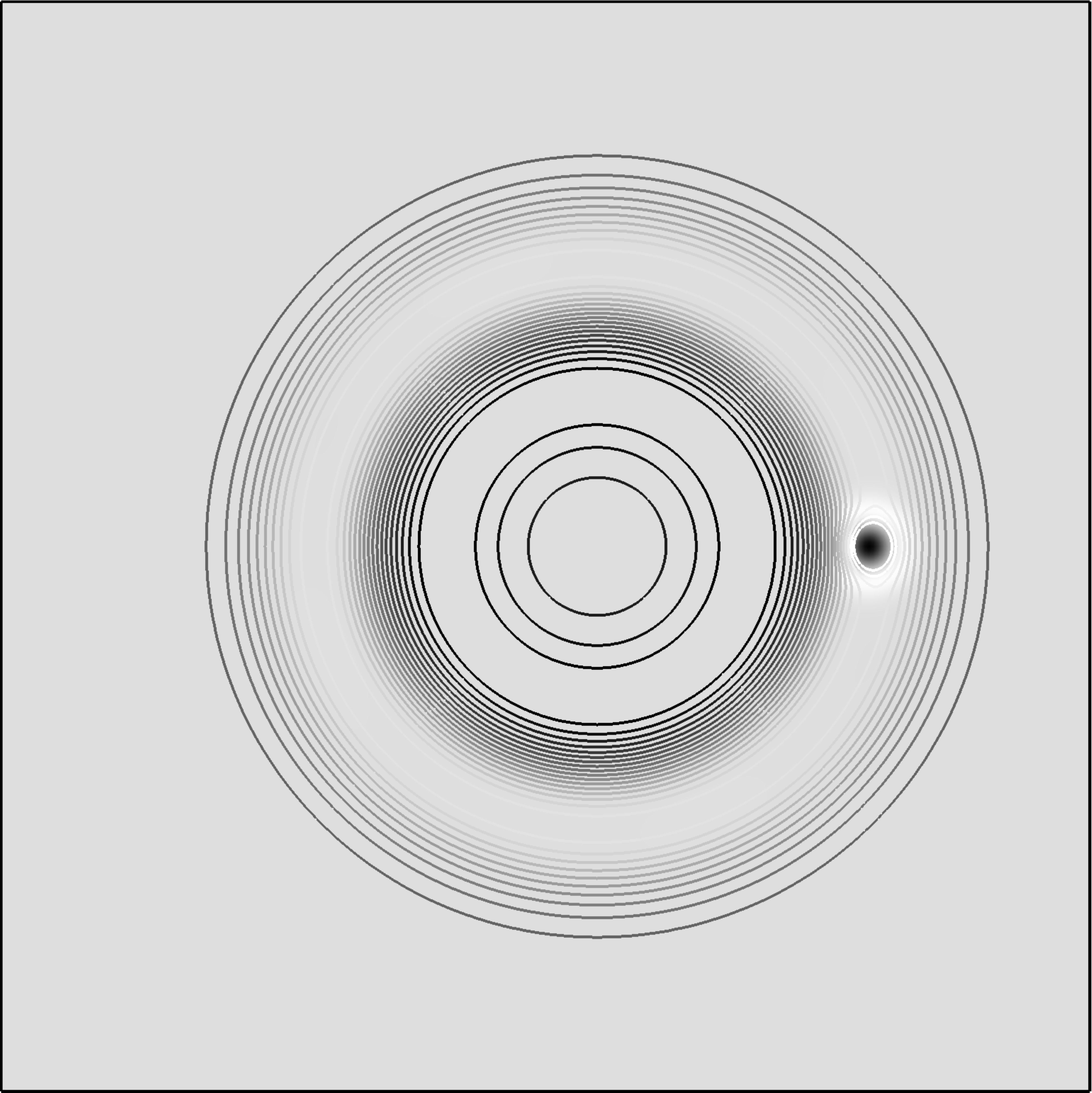}}
  \subfigure[Time $0.4$~s]{\includegraphics{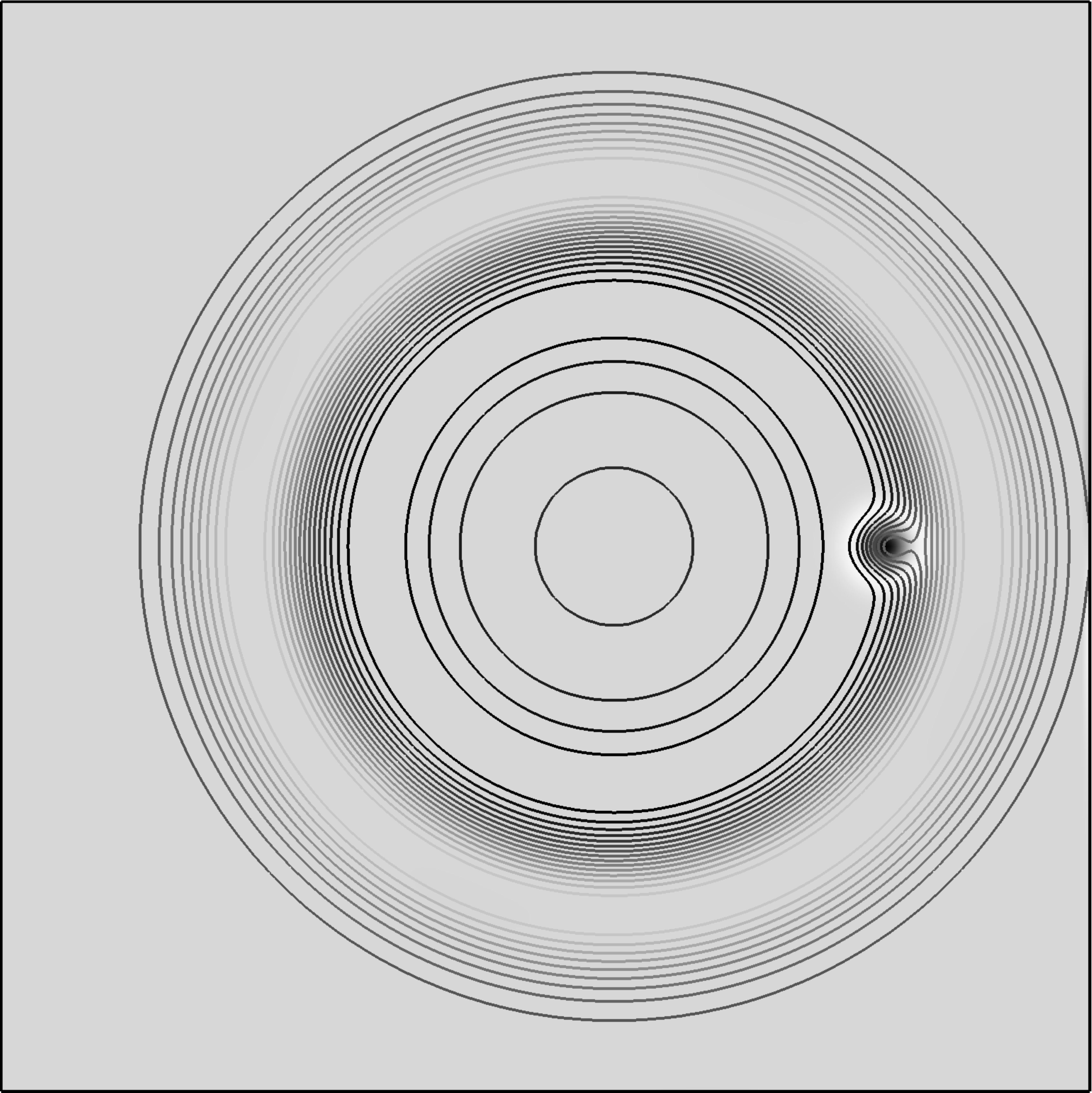}}
 \end{subfigmatrix}
 \caption{Isocontour of density superimposed on field of vorticity for solutions at $t=0.1$~s, $0.2$~s, $0.3$~s and $0.4$~s. The top row (figures (a)-(d)) are solutions computed with the purely low-Mach-number approach. The bottom row (figures (e)-(h)) are solutions computed with the hybrid method detailed in the present paper.}
 \label{fig:results_2D_test_case_density}
\end{figure}

Time evolution of the solution is presented in Figure~\ref{fig:results_2D_test_case_density}. 
Figure~\ref{fig:results_2D_test_case_density}(a)-(d) in the top row are the solutions computed with the purely low-Mach-number approach, whereas Figure~\ref{fig:results_2D_test_case_density}(e)-(h) 
are solutions computed with the new hybrid method.   The compressible solution gives results visually indistinguishable from the hybrid approach so those are not shown here.  In both the hybrid and compressible solutions, the circular pressure wave generated from the center of the domain propagates in all directions. As the sound speed is far higher than the mean flow velocity, the acoustic wave passes the entropy pulse and eventually leaves the domain at $0.4$~s. When the purely low-Mach-number approach is employed, the pressure pulse in the center of the domain is considered as an entropy pulse, and is convected in the same way as the entropy pulse localized downstream. It is noted that the hybrid solution correctly captures the behavior of the waves generated from acoustic pulse despite the fact that the compressible grid under the acoustic pulse is at lower resolution than in the fully compressible solution, and has an overset fine low-Mach-number grid.

In order to provide quantitative results, both the solution computed with the hybrid method and the purely compressible solution are compared to a reference exact analytical solution \cite{Tam:1993}. 
The numerical error is assessed by computation of the $\mathcal{L}^2$-norm of the difference between the computed and the reference solutions, which is expressed as follows:
\begin{equation}
\varepsilon_{\phi} = \mathcal{L}_{\phi}^2 \left(S_{sol} - S_{ref} \right) = \sqrt{\frac{\left(\phi_{sol} - \phi_{ref} \right)^2}{N_x N_y}}
\end{equation}
where subscripts $sol$ and $ref$ identify the computed and reference solutions, $\phi$ is the variable investigated, and $N_x$ and $N_y$ 
are the number of points in the $x$ and $y$ directions. Note that for simplicity, $N_x=N_y$.

Similarly to Sec.~\ref{subsec:results_1d_acoustic}, simulations are performed on a multi-levels grid set composed by a total of $L=5$ levels.  The first level $l=1$ is discretized with $N^{l=1}_x=32$ and $N^{l=1}_y=32$ points, while the other levels are progressively discretized with a mesh refinement ratio of a factor of $2$. Table~\ref{tab:summary_2D_waves_meshgrid} presents the configuration of the multi-levels grid set by providing a summary of $N_x$ and $N_y$ for each level $l$ of mesh refinement. 

\begin{table}
\centering
\renewcommand\arraystretch{1.5}
\begin{tabular}{c || c | c | c | c | c  }
\renewcommand*{\arraystretch}{0.5}
$l$    & $1$  & $2$ & $3$ & $4$ & $5$   \\ \hline\hline 
$N_x$  & $32$ & $64$ & $128$ & $256$ & $512$   \\\hline
$N_y$  & $32$ & $64$ & $128$ & $256$ & $512$    
\end{tabular}
\caption{Summary of the configuration for simulations performed on the 2D mixed modes propagation test case.}
\label{tab:summary_2D_waves_meshgrid}
\end{table}

Simulations are performed by first selecting, from $l=1$ to $l=4$, the level $l_{\rm Comp}$ where the fully compressible Eqs.~(\ref{eqn:mass_comp}-\ref{eqn:energy_comp}) are solved, and then by selecting a successive addition of low-Mach-number levels of mesh refinement, the finest level being designed by $L$. In total, $10$ simulations are performed, and the choices of $l_{\rm Comp}$ and $L$ for each simulation are summarized in Table~\ref{tab:summary_2D_waves}. Furthermore, the time-steps for both the compressible and low-Mach-number equations are computed as described in Sec.~\ref{subsubsec:step1} and $\varepsilon_{\phi}$ is computed for solutions taken at the time $t=0.3$~s.

\begin{table}
\centering
\renewcommand\arraystretch{1.5}
\begin{tabular}{c || c | c | c | c | c  }
\renewcommand*{\arraystretch}{0.5}
\diagbox{$l_{\rm Comp}$}{$L$} & $1$ & $2$ & $3$ & $4$ & $5$   \\ \hline\hline 
$1$ &  & $\times$ & $\times$ & $\times$ &  $\times$  \\\hline
$2$ &  &   & $\times$ & $\times$ &  $\times$  \\\hline
$3$ &  &  &  & $\times$ &  $\times$  \\\hline
$4$ &  &  &  &  &  $\times$  \\
\end{tabular}
\caption{Summary of the choices of $l_{\rm Comp}$ and $L$ for all simulations performed during spatial convergence test of the hybrid method with the propagation of mixed acoustic, entropic and vorticity modes in a 2D square domain.}
\label{tab:summary_2D_waves}
\end{table}

\begin{figure}[!ht]
 \begin{subfigmatrix}{2}
  \subfigure[]{\includegraphics{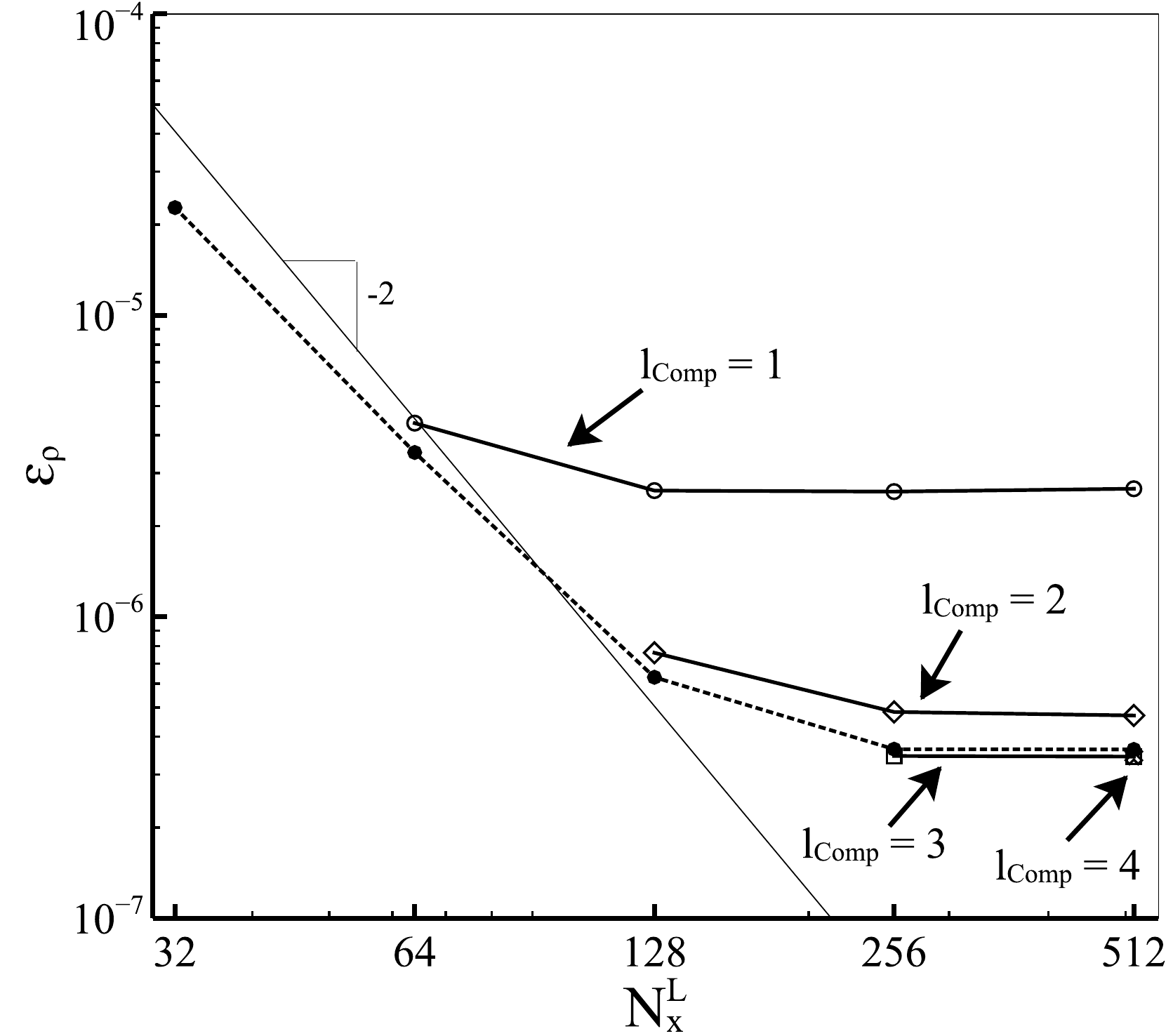}}
  \subfigure[]{\includegraphics{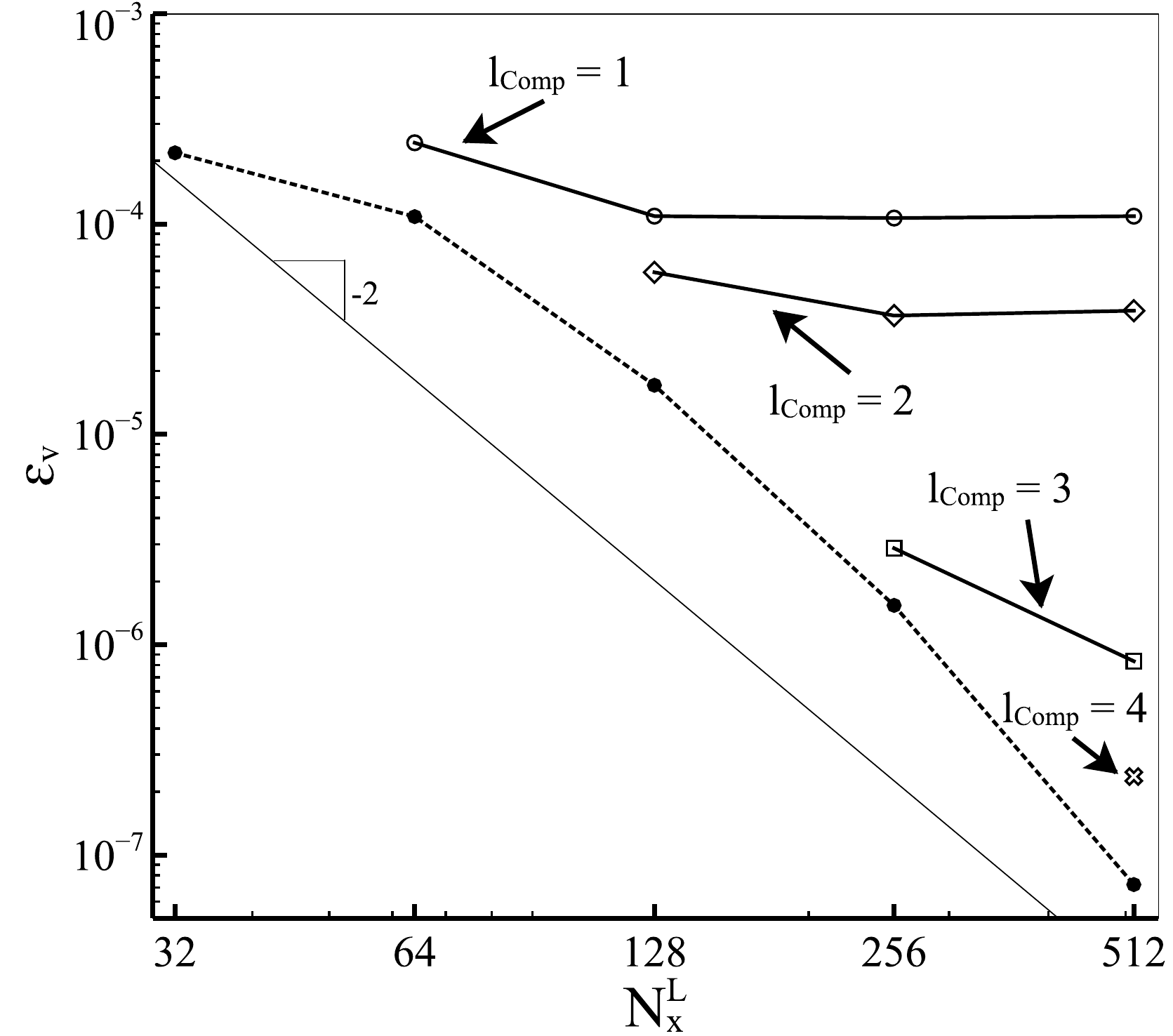}}
  \end{subfigmatrix}
 \caption{$\mathcal{L}^2$-norm of the discretization error for different maximum level $L$ of mesh refinement for the low-Mach-number equations: (a) $\varepsilon_{\rho}$ for the density, (b) $\varepsilon_{v}$ for the velocity in the $y$-direction. Circle, diamond, square and cross symbols represent the fully compressible equations solved on the level $l_{\rm Comp}$ set at $l=1$, $l=2$, $l=3$ and $l=4$, respectively. The dashed black line represents the evaluation of $\varepsilon_{\rho}$ and $\varepsilon_{v}$ for the purely compressible approach. The solid black line represents a second order slope.}
 \label{fig:results_2D_acoustic_wave_sptial_convergence_rate}
\end{figure}

Figures~\ref{fig:results_2D_acoustic_wave_sptial_convergence_rate}.(a) and~\ref{fig:results_2D_acoustic_wave_sptial_convergence_rate}.(b) present the $\mathcal{L}^2$ norm error computed for the density $\left(\varepsilon_{\rho}\right)$ and the velocity in the $y$-direction $\left(\varepsilon_{v}\right)$, respectively. Circle, diamond, square and cross symbols represent $l_{\rm Comp}$ set at $l=1$, $l=2$, $l=3$ and $l=4$, respectively. This corresponds to a discretization of $N_x^{l_{\rm Comp}}=32, 64, 128$ and $256$ points, respectively.  Moreover, the dashed lines represent $\varepsilon_{\rho}$ and $\varepsilon_{v}$ evaluated from the solutions computed with the purely compressible equations, while the solid line is the second order slope.  

Note here that $\varepsilon_{\rho}$ and $\varepsilon_{v}$ are not computed in the full 2D domain but only on the $x-$axis taken at $y=L_{x}/2$. This specific choice enable us to separate the contribution of acoustic, entropic and vorticity modes. Indeed, as the axis is taken along the propagation of the acoustic wave, no contribution from the acoustic and entropic modes should appear in the $v$ component of the velocity, but only the ones from the vortex structure. In contrary, on this specific axis, only acoustic and entropic modes should contribute to the evaluation of the density, and not the vorticity mode. 

In Figure.~\ref{fig:results_2D_acoustic_wave_sptial_convergence_rate}.(a), the evaluation of $\varepsilon_{\rho}$ for the solutions computed with the purely compressible equations (dashed line) follows a second order rate of convergence, and starts to reach a plateau for levels $l > 3$ (viz. $N_x^L > 128$). When the hybrid method is employed, the contribution of solving the low-Mach-number equations on an additional level significantly reduces $\varepsilon_{\rho}$ to approximately get the same error as if the additional layer was employed to solve the fully compressible equations. However, solving the low-Mach-number equations on additional finest levels does not help significantly to further reduces $\varepsilon_{\rho}$, which also reach eventually a plateau. This suggest that  solving the low-Mach-number equations on additional levels of mesh refinement strongly reduced the error made on the convection of the entropy spot, but that the numerical errors made because of the poor resolution of the acoustic wave on the coarser mesh still remain in the solution at the finest level. This statement is in accordance with the convergence rate behavior observed in Sec.~\ref{subsec:results_1d_acoustic} for the propagation of purely acoustic waves.

Furthermore, the same observations can be made from Figure~\ref{fig:results_2D_acoustic_wave_sptial_convergence_rate}.(b). Recall that only contributions from the vorticity mode should appear in the solution, solving the low-Mach-number equations on additional finer levels should strongly reduce $\varepsilon_{v}$. However, a significant error remains on $\varepsilon_{v}$ when $l_{\rm Comp}=1$ and $2$, even at the finest level of refinement for the low-Mach-number equations. This suggests that numerical errors from the poor resolution of the acoustic wave appear in the low-Mach-number solution. For $l_{\rm Comp}=3$, the acoustic wave is considered enough well resolved, so that numerical errors from the purely compressible equations become negligible and the contribution of additional low-Mach-number levels is significant to reduce the overall error made on the velocity. This is consistent with the observation made in Figure~\ref{fig:results_2D_acoustic_wave_sptial_convergence_rate}.(a) that the error in the density has reached a plateau for $l_{\rm Comp}>3$. 

As a partial conclusion, this study exhibits the limitations of the hybrid method. Solving the low-Mach-number equations on additional level of mesh refinement only provides a better solution for phenomena that do not include contributions from the acoustics. This suggests that acoustic phenomena of interest must still be well enough resolved on the levels where the purely compressible equations are solved. This is obvious with the present test case. For example in Figure~\ref{fig:results_2D_acoustic_wave_sptial_convergence_rate}.(a), for $l_{\rm Comp}=3$ and $4$, the hybrid method provides an error $\varepsilon_{\rho}$ that is similar to the error made with the purely compressible approach (dashed line). 

\begin{figure}[!ht]
\centering
  \includegraphics[width=0.65\textwidth]{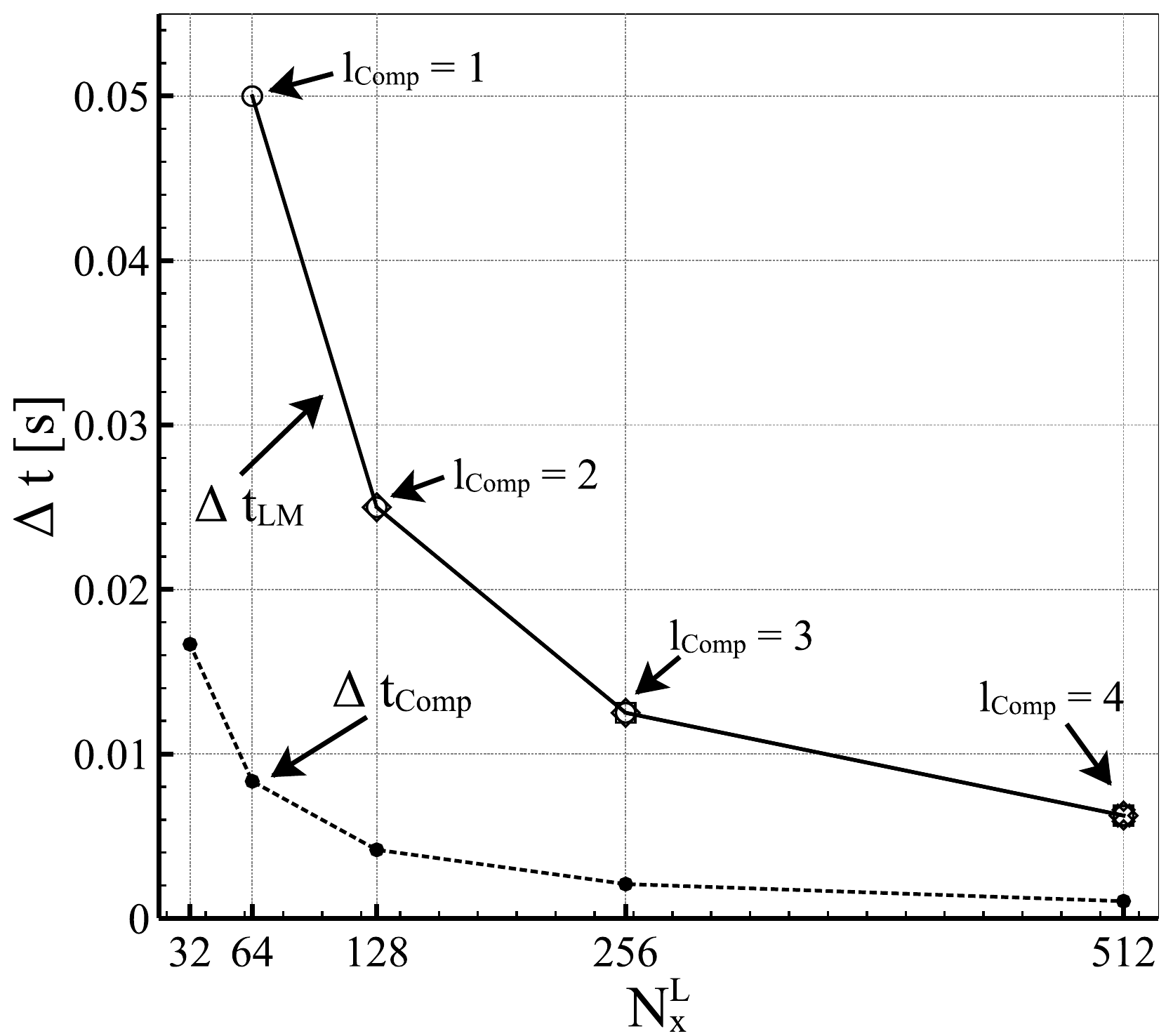}
 \caption{Average time-step employed during simulations performed with the purely compressible approach (dashed line) and the hybrid method (symbols), and for different maximum level $L$ of mesh refinement for the low-Mach-number equations. For the hybrid method, circle, diamond, square and cross symbols represent $l_{\rm Comp}$ set at $l=1$, $l=2$, $l=3$ and $l=4$, respectively.}
 \label{fig:results_2D_waves_dt_vs_Nx}
\end{figure}

\begin{figure}[!ht]
\centering
\includegraphics[width=0.65\textwidth]{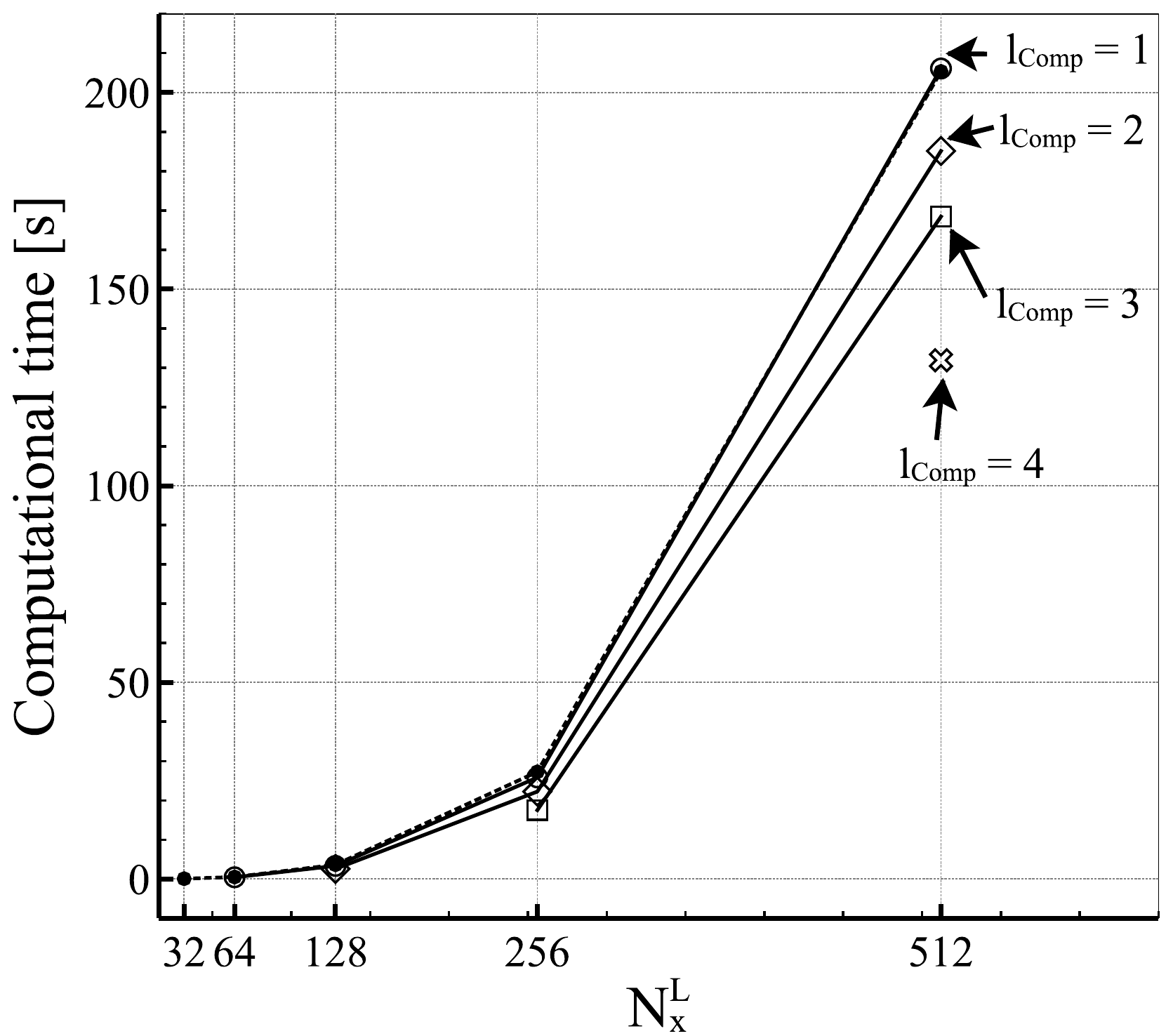}
 \caption{Wall-clock computational time spent to perform simulations with the purely compressible approach (dashed line) and the hybrid method (symbols), and for different maximum level $L$ of mesh refinement for the low-Mach-number equations. For the hybrid method, circle, diamond, square and cross symbols represent $l_{\rm Comp}$ set at $l=1$, $l=2$, $l=3$ and $l=4$, respectively.}
 \label{fig:results_2D_waves_comput_time_vs_Nx}
\end{figure}

However, the interest of the hybrid method developed in the present paper is highlighted in Figures~\ref{fig:results_2D_waves_dt_vs_Nx} and \ref{fig:results_2D_waves_comput_time_vs_Nx}. Figure~\ref{fig:results_2D_waves_dt_vs_Nx} presents the comparison of the average time-step employed during simulations performed with the purely compressible approach (dashed line) and the hybrid method (symbols). For the hybrid method, similarly to  Figures~\ref{fig:results_2D_acoustic_wave_sptial_convergence_rate}.(a) and~\ref{fig:results_2D_acoustic_wave_sptial_convergence_rate}.(b), the circle, diamond, square and cross symbols represent $l_{\rm Comp}$ set at $l=1$, $l=2$, $l=3$ and $l=4$, respectively. They obviously collapse in the same curve because the finest level of mesh refinement $L$ determines the low-Mach-number time-step $\Delta t_{\rm LM}$. On the other hand, Figure~\ref{fig:results_2D_waves_comput_time_vs_Nx} presents the overall wall-clock computational time corresponding to the simulations performed in the present section. Together with the results presented in Figure~\ref{fig:results_2D_waves_dt_vs_Nx} and Figures~\ref{fig:results_2D_acoustic_wave_sptial_convergence_rate}.(a) and ~\ref{fig:results_2D_acoustic_wave_sptial_convergence_rate}.(b), two major general observations can be made:

\begin{itemize}
\item When $l_{\rm Comp}$ is too coarse, solving the low-Mach-number equations on additional levels of mesh refinement does not help to capture a good representation of the physics, or to provide a significant gain in the computational time.
\item once the physics specifically related to generation of the acoustics is well enough resolved by selecting the proper level of discretization $l_{\rm Comp}$, solving the low-Mach-number equations on a few additional levels provides a significant gain in the computational effort, while providing lower numerical errors in the solution. This is particularly true for the configuration $l_{\rm Comp}=4$ and $L=5$: the hybrid method provides a discretization error in the density which is lower than the purely compressible approach, while at the same time exhibiting a computational cost about twice less expensive. Note that the reduction in numerical errors is strongly dependent of the problem simulated, as well as the procedure employed for adaptive discretization of the flow.
\end{itemize}

Note that in Figure~\ref{fig:results_2D_waves_dt_vs_Nx}, the time-steps employed by the hybrid method are significantly larger than the ones computed by the fully compressible approach. However, in Figure~\ref{fig:results_2D_waves_comput_time_vs_Nx}, one can observe that the gain in the computational time provided by the hybrid method becomes significant for $l_{\rm Comp}>3$. This can be explained by the fact that, as the tolerance parameter $\epsilon_p$ in Eq.~(\ref{eqn:epsilon_p}) is set to $\epsilon_p = 1 \times 10^{-12}$, many sub-iterations are required (approximately $m=20$) when $l_{\rm Comp}$ is too coarse, because the fine low-Mach-number solution deviates significantly from the badly resolved compressible solution. However when the acoustics is well resolved enough, for example for $l_{\rm Comp}=3$, it has been observed that the low-Mach-number solution converges very quickly to the compressible solution, in a few iterations (on average, approximately $m=2$).

The present test case highlights the capacity of the hybrid method to retain acoustic phenomena within the context of a low-Mach-number solver. The major trend highlighted in this section is that acoustic phenomena must be well enough resolved where the fully compressible equations are solved. It is however emphasized that this test case is very canonical because the acoustics and the rest of the dynamic of the flow are, in the same time, well defined and decoupled from each other. For practical applications, the goal is to solve the low-Mach-number equations only in regions of the domain where the Mach number is small -- hence the computational savings due to the larger low-Mach-number time step are greatest -- and where the flow features have very fine structure that must be resolved. This practical application is now investigated in the following section by the computation of the aeroacoustic sound generated by the vortex formation from a Kelvin-Helmholtz instability in low-Mach-number mixing layers.

\subsection{Aeroacoustic propagation from a low-Mach-number Kelvin-Helmholtz instability}
\label{subsec:results_2d_mixing_layers}

The present test case aims to evaluate the performance of the hybrid method developed in this paper for a realistic physical phenomenon that can appear in practical flow applications similar to the ones encountered in the industry. A Kelvin-Helmholtz instability in low-Mach-number mixing layers is simulated. Basically, the interface between two flows in opposite directions is excited on the most unstable mode of fluctuations. A series of small vorticity structures  progressively appear, before eventually merging into a single rotating vortex. As vortex breaking is a source of aeroacoustic sound, pressure waves are generated and propagate inside the domain. The key particularity of the present configuration is that the acoustic wavelength is large,   with a typical size of the order of half of a meter. In contrary, the mixing layer interface is very small, or the order of a millimeter. Consequently, there is a large disparity between the spatial scales of the vorticity structures and the aeroacoustic waves propagated in the domain.

While being a canonical test case with a well-controlled physics of the flow, this test case is representative of the phenomena that appear in the context of noise generated by jets in practical industrial applications. Therefore, this test case has been widely computed in the aeroacoustic community to understand the sources of vortex sound generation, as well as to evaluate the performances of computational aeroacoustic techniques as mentioned in the introduction part of the present paper (see \cite{Bogey:2002,Fortune:2004,Golanski:2004,Golanski:2005}, among others). Indeed, the main issue here is that the mixing interface must be well enough resolved in order to capture accurately the vortex formation, which is critical to capture as well the proper aeroacoustic phenomena, especially in terms of frequency and pressure amplitudes. Consequently, this test case is a good candidate to assess the performance of the hybrid method developed in the present paper.

 The configuration of the test case is inspired by the temporal representation of the instability as proposed by Golanski \emph{et al.}\cite{Golanski:2004}, which features a controlled excitation to generate several pairs of vortices that eventually merge together and generate noise. The computational domain is a rectangle of dimension $L_x \times L_y$, with $L_x = 2 \lambda_a$ and $L_y=64 \lambda_a$. Here, according to the linear stability theory \cite{Michalke:1964,Sandham:1991} , $\lambda_a=\frac{2 \pi}{k_a}\delta_{\omega}$ is the wavelength of the most unstable mode in the mixing interface, where $k_a=0.4446$ is the wavenumber of maximum amplification and $\delta_{\omega}$ is the thickness of the mixing layers interface. The initial flow conditions are given as follows:

\begin{align}
\rho^{\rm init}\left(x,y\right) &= \rho_{\rm ref} \label{eqn:2d_vortex_1}  \\
u^{\rm init}_x\left(x,y\right) &= \frac{U_1 + U_2}{2} + \frac{U_1 - U_2}{2} \tanh\left(\frac{2 \left(y - y_{\rm ref} \right) }{\delta_{\omega}}\right) \\ u^{\rm init}_y\left(x,y\right) &= A e^{- \sigma \left(\frac{y-y_{\rm ref}}{\delta_{\omega}} \right)^2} \times \left[\cos \left(\frac{8 \pi}{L_x} x\right) + \frac{1}{8}\cos \left(\frac{4 \pi}{L_x} x\right)+\frac{1}{16}\cos \left(\frac{2 \pi}{L_x} x\right) \right] \\
p_0^{\rm init}\left(x,y\right) &= p_{\rm ref}, \hspace{0.5cm} p_1^{\rm init}\left(x,y\right) = 0 \label{eqn:2d_vortex_2}
\end{align}
Here, $\rho_{\rm ref} = 1.1$~kg/m$^3$ and $p_{\rm ref}=9\times10^5$~Pa, while $\gamma=1.1$ so that the speed of sound is $c_{\rm ref}=300$~m/s. The mean velocity of the lower and upper flows are set to $U_1 = 20$~m/s and $U_2 = - U_1$, respectively. The thickness of the mixing layers interface is defined by $\delta_{\omega}=1 \times 10^{-3}$~m. The parameters $A=0.025 \left(U_1 - U_2\right)$ and $\sigma=0.05$ control the amplitude and the thickness of the perturbation imposed to the mean flow field. Finally, $y_{\rm ref}=L_y/2$ is set so as to center the mixing layers interface in the middle of the domain. Overall, the mean Mach number in the simulation is approximately $M \approx 0.06$. Note that in order to impose a divergence-free initial condition, a projection in pressure is initially performed. Basically this operation is similar to solving Eq.~(\ref{eqn:poisson_eq}) and~(\ref{eqn:advance_u_np1}), but with $\nabla \mathbf{u}^{n+1}=0$ and $\mathbf{u}^{*}$ being the initial flow provided by Eqs.~(\ref{eqn:2d_vortex_1})-(\ref{eqn:2d_vortex_2}). Finally, the tolerance parameter $\epsilon_p$ in Eq.~(\ref{eqn:epsilon_p}) is set to $\epsilon_p = 1 \times 10^{-10}$, which corresponds to an average number of sub-iterations $m=5$.

\begin{figure}[!ht]
\centering
\includegraphics[width=0.65\textwidth]{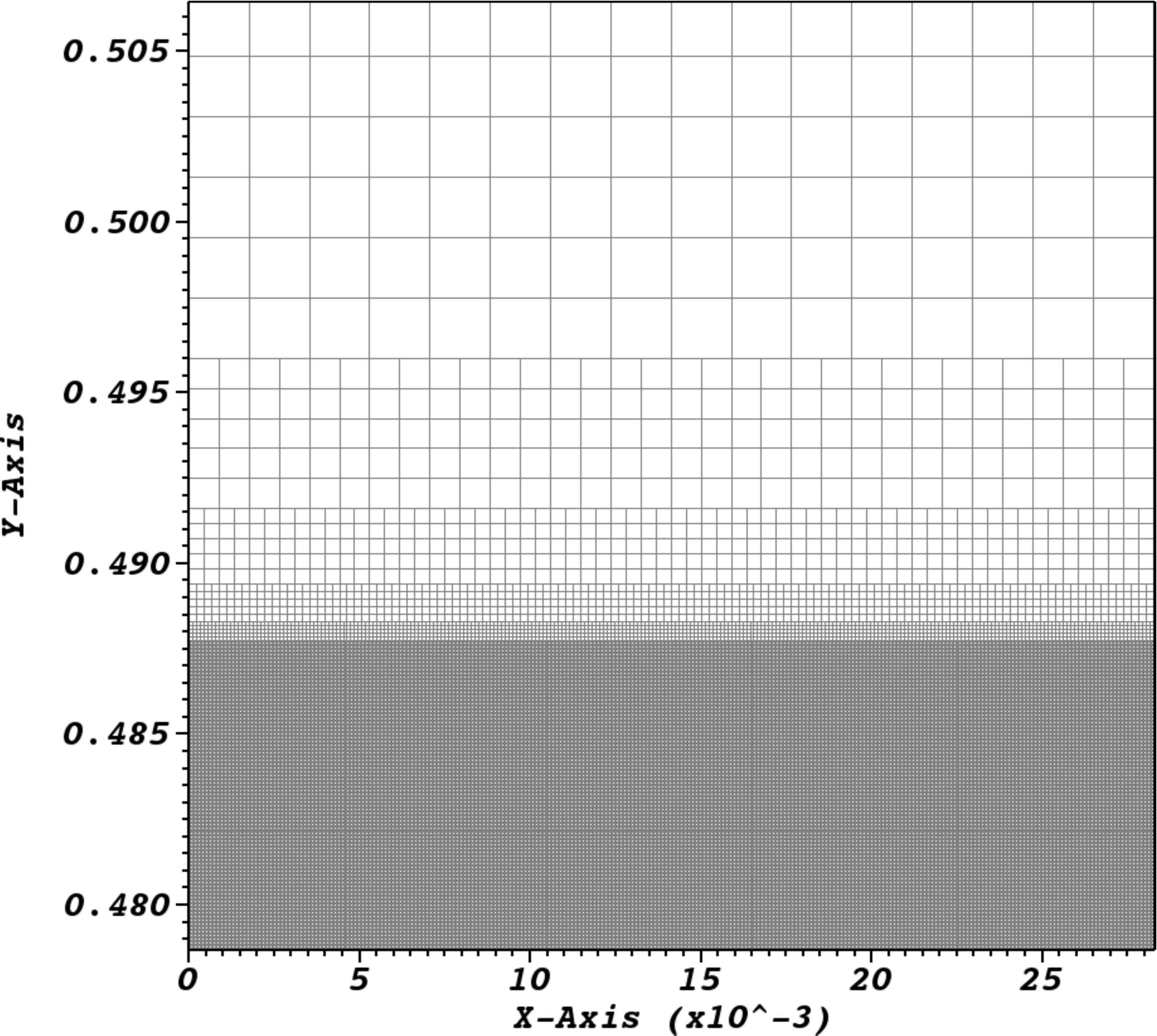}
 \caption{Representation of the multi-levels grid set around $L_y = 28 \lambda_a$.}
 \label{fig:results_2D_aeroacoustics_mesh_grid}
\end{figure}

In order to demonstrate the performances of the hybrid method developed in the present paper, the low-Mach-number Kelvin-Helmholtz instability case is simulated first with the purely compressible approach. As explained before, the mixing layers interface must be well enough resolved to accurately capture the vortex formation, but the acoustic waves exhibit a long wavelength that does not require such a fine discretization. In order to save computational resources, the Adaptive Mesh Refinement (AMR) framework is adopted. Note that here, for simplicity, the  additional mesh levels of refinement are imposed manually in the simulation, but they could have been specified by a criterion based on the vorticity for example.  Let us define $l_{\rm max\_comp}$  the total number of levels of mesh refinement. The whole domain is covered by a first level $l_{\rm max\_comp}=1$ consisting of very coarse grid,  defined as $N_x^{l=1}=16$ and $N_y^{l=1}=512$. This corresponds to a spatial grid size of $\Delta x = 1.76 \times 10^{-3}$~m. Recall that $\delta_{\omega}=1.0 \times 10^{-3}$~m, the mixing layers interface is then represented by barely 2 points, which is obviously too coarse to capture the vortex formation. 
Additional levels with a refinement factor of $2$ are successively superimposed on top of each other in the area of the computational domain comprised between $L_y = 28 \lambda_a$ and $L_y = 36 \lambda_a$. This area is selected so as to cover the full vortex evolution. As shown later, a total of $7$ additional levels of mesh refinement are required to capture accurately the formation of the vortex and to provide converged results in term of pressure evolution. The multi-levels grid set is depicted in Figure~\ref{fig:results_2D_aeroacoustics_mesh_grid}. Note that for each level of mesh refinement, a buffer zone of $4$ cells is imposed so as to let the solution to adapt between each level. Moreover, Table~\ref{tab:summary_2D_mixing_interface} summarizes, for each $l_{\rm max\_comp}$, the corresponding minimum $\Delta x$ and an approximation of the associated numbers $N_{\rm interface}$ of grid points in the mixing layers interface.

\begin{table}
\centering
\renewcommand\arraystretch{1.5}
\begin{tabular}{c || c | c   }
\renewcommand*{\arraystretch}{0.5}
$l_{\rm max\_comp}$ & $\Delta x$ in $[m]$ & $N_{\rm interface}$   \\ \hline\hline 
$1$ &  $1.765 \times 10^{-3}$ &  $2$ \\ \hline
$2$ &  $8.825 \times 10^{-4}$ &  $3$ \\ \hline
$3$ &  $4.412 \times 10^{-4}$ &  $5$ \\ \hline
$4$ &  $2.206 \times 10^{-4}$ &  $10$ \\ \hline
$5$ &  $1.103 \times 10^{-4}$ &  $20$ \\ \hline
$6$ &  $5.516 \times 10^{-4}$ &  $40$ \\ \hline
$7$ &  $2.758 \times 10^{-5}$ &  $80$ \\ \hline
$8$ &  $1.380 \times 10^{-5}$ &  $160$ \\ \hline
\end{tabular}
\caption{Summary, for each $l_{\rm max\_comp}$, of the corresponding minimum $\Delta x$ and an approximation of the associated numbers $N_{\rm interface}$ of grid points in the mixing layers interface.}
\label{tab:summary_2D_mixing_interface}
\end{table}

\begin{figure}[ht]
 \begin{subfigmatrix}{2}
  \subfigure[Time $0.5$~ms]{\includegraphics{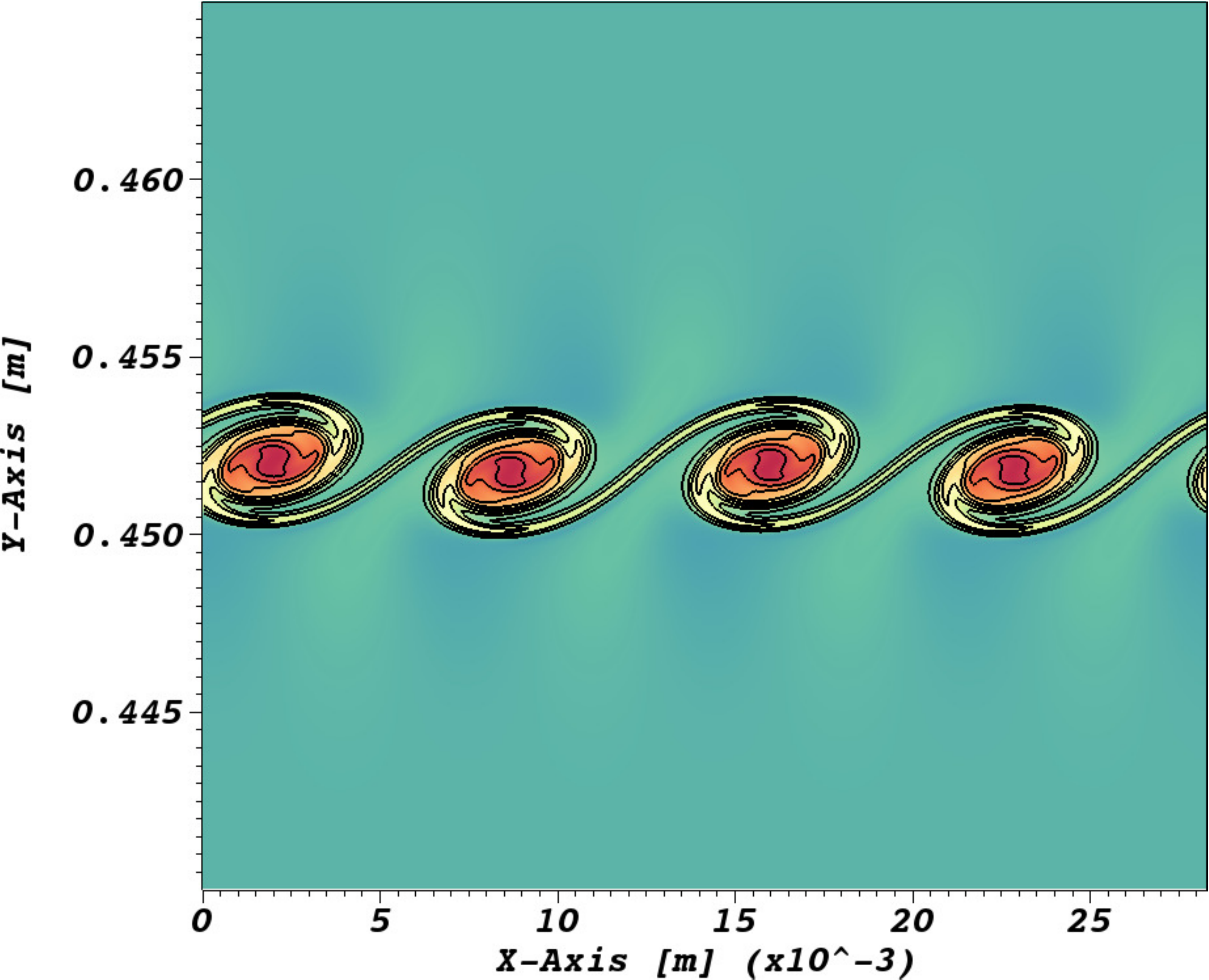}}
  \subfigure[Time $1.0$~ms]{\includegraphics{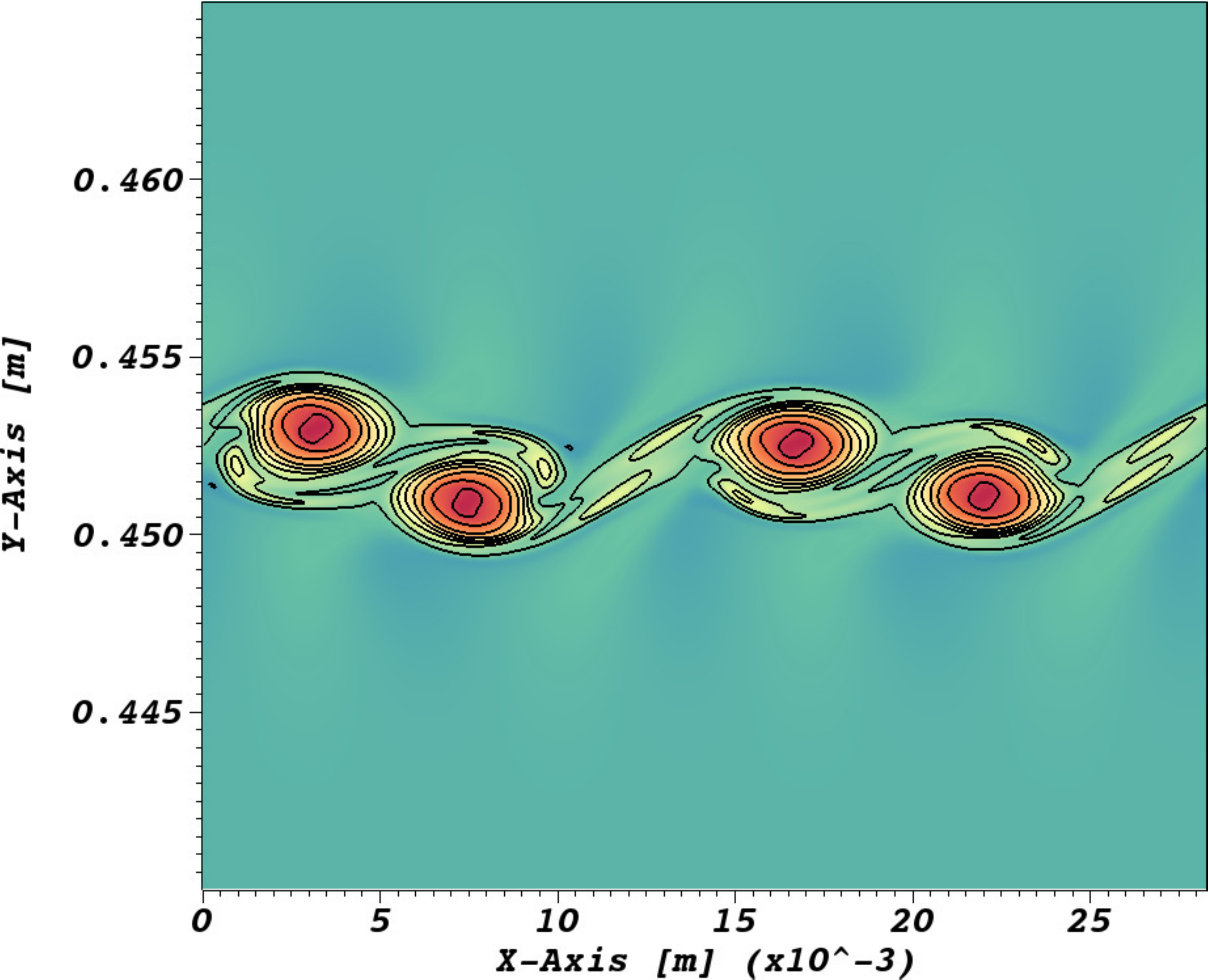}}
  \subfigure[Time $2.0$~ms]{\includegraphics{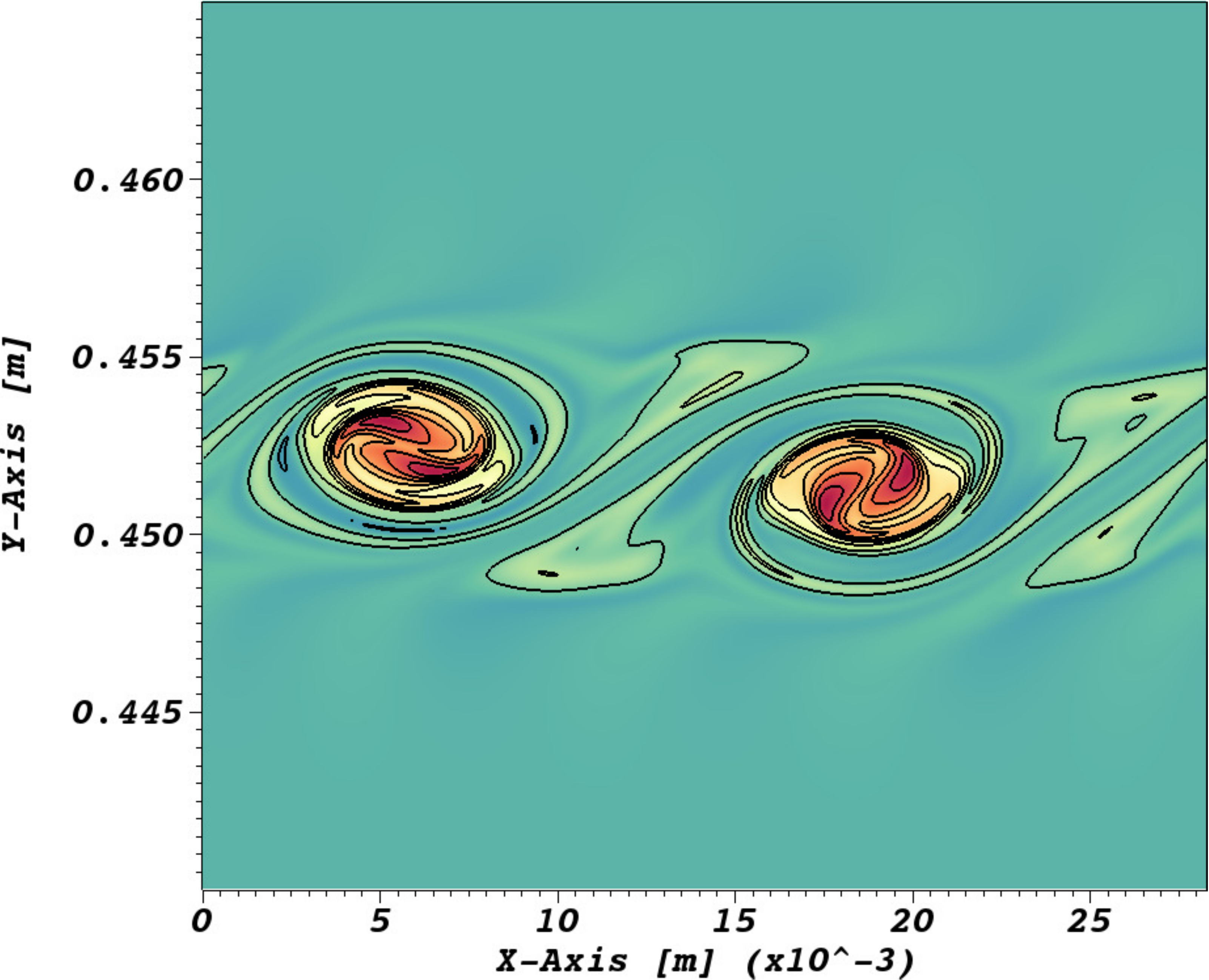}}
  \subfigure[Time $4.0$~ms]{\includegraphics{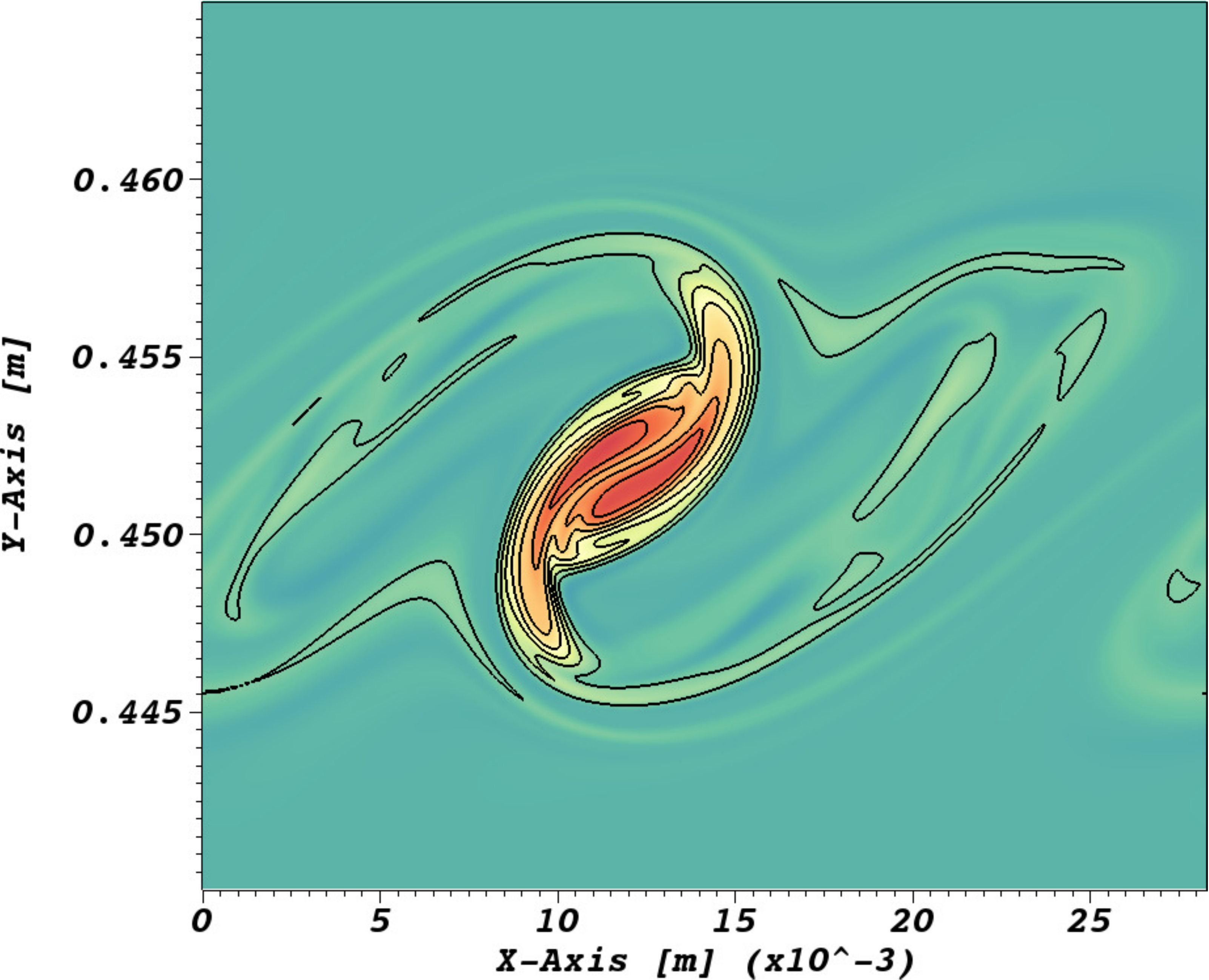}}
 \end{subfigmatrix}
 \caption{Fields of vorticity at different time of the simulation, computed with the purely compressible approach with $8$ levels of mesh refinement. Contours of vorticity are also depicted to visually identify the evolution of the mixing layers interface.}
 \label{fig:results_2D_aeroacoustics_vorticity}
\end{figure}

Simulations are performed over a time of $4 \times 10^{-3}$~s. Contours of the vorticity are depicted in Figure.~\ref{fig:results_2D_aeroacoustics_vorticity} for a selection of temporal snapshots. At $t=0.5 \times 10^{-3}$~s (see Figure.~\ref{fig:results_2D_aeroacoustics_vorticity}.(a)), the interface is still clearly visible but is distorted to form $4$ vortex structures. Very quickly, at  $t=1.0 \times 10^{-3}$~s (see Figure.~\ref{fig:results_2D_aeroacoustics_vorticity}.(b)), the vortex structures are merging together two by two (see Figure.~\ref{fig:results_2D_aeroacoustics_vorticity}.(c)), and these two structures then merge in a final unique rotating vortex (see Figure.~\ref{fig:results_2D_aeroacoustics_vorticity}.(d)). During this process, acoustic pressure is generated and propagates in the domain.

\begin{figure}[!ht]
\centering
\includegraphics[width=0.65\textwidth]{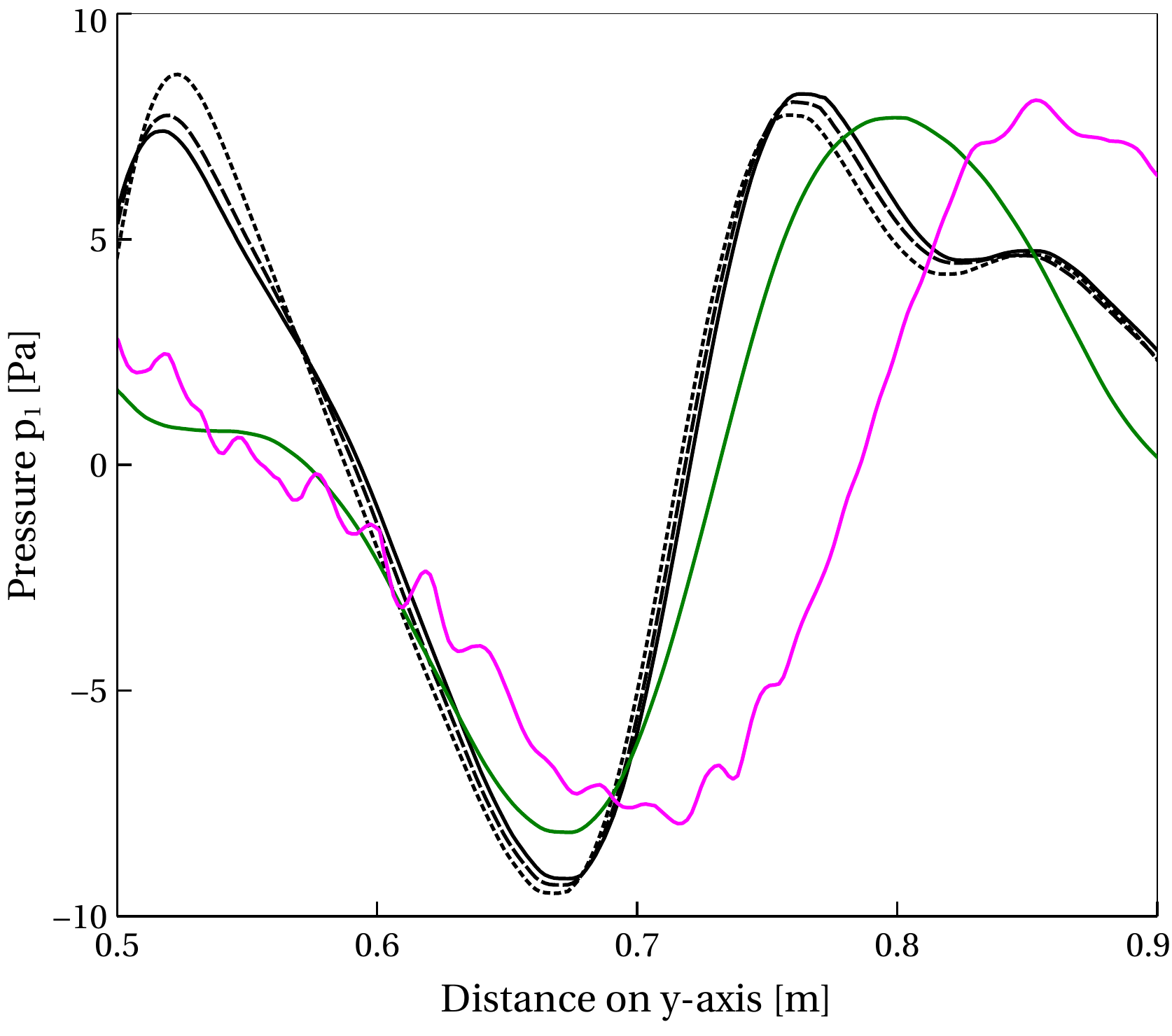}
 \caption{Signal of pressure fluctuations $p_1$ at $t=4 \times 10^{-3}$~s taken on the $y-$axis in the upper part of the domain between $L_y = 36 \lambda_a$ and $L_y = 64 \lambda_a$. Solutions computed with the purely compressible approach with  $l_{\rm max\_comp}=2$ (magenta solid line), $l_{\rm max\_comp}=4$ (green solid line), $l_{\rm max\_comp}=6$ (dashed black line), $l_{\rm max\_comp}=7$ (dotted black line) and $l_{\rm max\_comp}=8$ (solid black line).}
 \label{fig:results_2D_aeroacoustics_pressure_pure_comp}
\end{figure}

Figure~\ref{fig:results_2D_aeroacoustics_pressure_pure_comp} presents the signal of pressure fluctuations $p_1$ at $t=4 \times 10^{-3}$~s taken on the y-axis in the upper part of the domain, namely between $L_y = 36 \lambda_a$ and $L_y = 64 \lambda_a$, and for different levels of mesh refinement. The solid magenta line in Figure~\ref{fig:results_2D_aeroacoustics_pressure_pure_comp} represents the pressure for $l_{\rm max\_comp}=2$. As reported in Table~\ref{tab:summary_2D_mixing_interface}, this correspond to a spatial grid size in the mixing layers interface is $\Delta x = 8.825 \times 10^{-4}$~m, i.e approximatively $3$ points in the mixing layer. The green solid line represents the pressure for $l_{\rm max\_comp}=4$, while the black dotted and dashed lines corresponds to $l_{\rm max\_comp}=6$ and $l_{\rm max\_comp}=7$, respectively. Finally, the solid black line corresponds to $l_{\rm max\_comp}=8$ and is considered as a converged solution.  This corresponds to distribution of $160$ points in the initial mixing layers interface thickness. It is quite obvious here that a coarse discretization of the interface leads to a very poor representation of the acoustic wave, especially in terms of the associated frequency and phase relationship with the vortex.

\begin{figure}[!ht]
\centering
\includegraphics[width=0.65\textwidth]{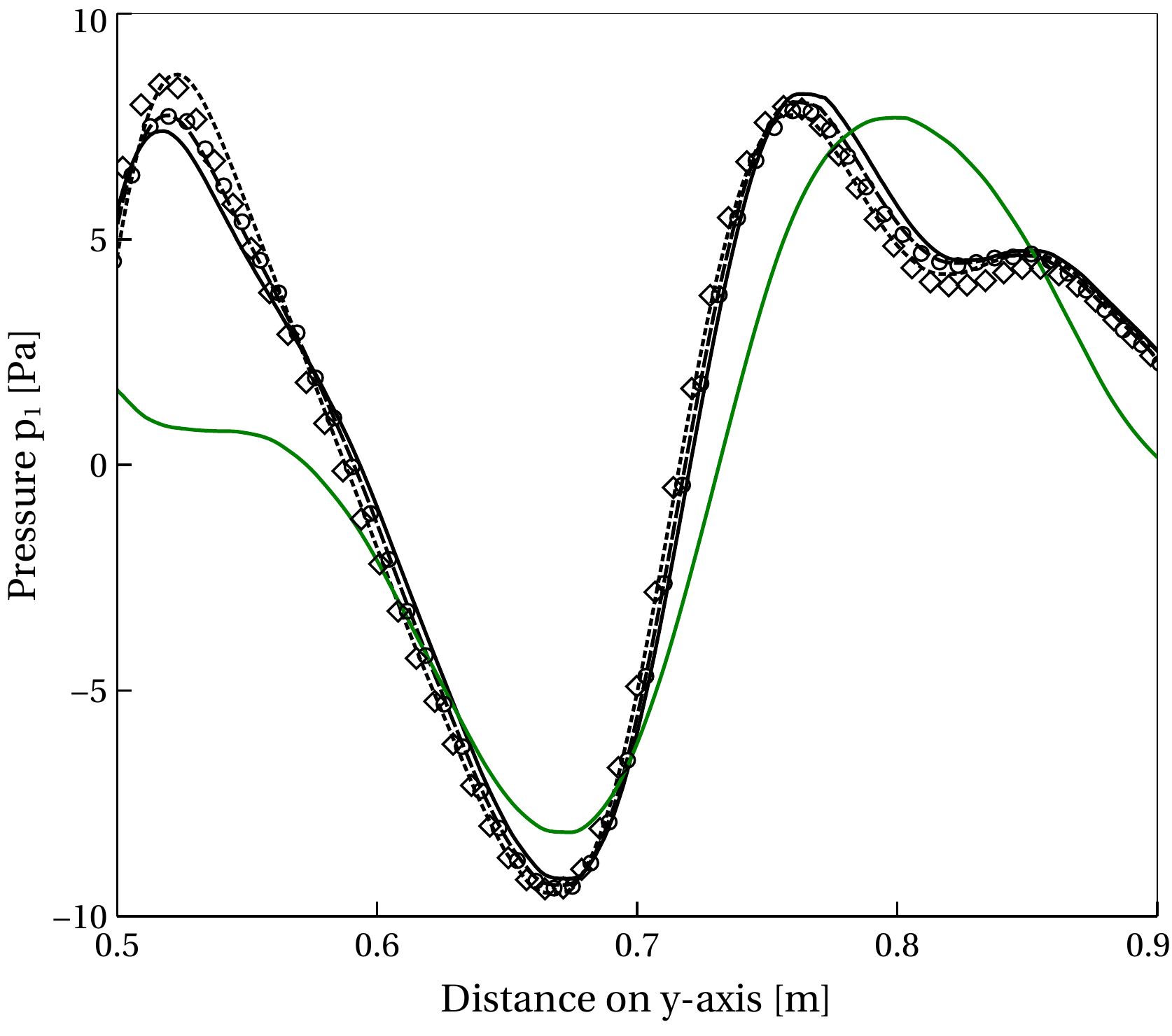}
 \caption{Signal of pressure fluctuations $p_1$ at $t=4 \times 10^{-3}$~s taken on the $y-$axis in the upper part of the domain between $L_y = 36 \lambda_a$ and $L_y = 64 \lambda_a$. Solutions computed with the purely compressible approach with   $l_{\rm max\_comp}=4$ (green solid line), $l_{\rm max\_comp}=6$ (dashed black line), $l_{\rm max\_comp}=7$ (dotted black line) and $l_{\rm max\_comp}=8$ (solid black line). Solutions computed with the hybrid method with $l_{\rm max\_comp}=4$ and $L=5$ (square symbols) and $L=6$ (circle symbols).}
 \label{fig:results_2D_aeroacoustics_pressure_hybrid2}
\end{figure}

The present configuration is now simulated with the hybrid method described in this paper. Again, the signal of pressure fluctuations $p_1$ at $t=4 \times 10^{-3}$~s is taken on the $y-$axis in the upper part of the domain. Results are gathered in Figure~\ref{fig:results_2D_aeroacoustics_pressure_hybrid2}. The colors and shapes of the lines are the same as in Figure~\ref{fig:results_2D_aeroacoustics_pressure_pure_comp} and corresponds to the results with the purely compressible approach. The symbols correspond to the results computed with the hybrid method. For all simulations performed with the hybrid method, $l_{\rm max\_comp}=4$. The square and circle symbols correspond to the results when the low-Mach-number equations are solved on $1$ and $2$ additional layers of mesh refinement, respectively. Quantitative results are presented in Table~\ref{tab:summary_2D_L2_pressure}. The left column the $\mathcal{L}_2$-norm of the error $\varepsilon_{p}$ computed at $t=4 \times 10^{-3}$~s for the pressure $p_1$ between simulations performed either with the hybrid method or the fully compressible approach at different levels $l_{\rm max\_comp}=1, \ldots, 7$, and the reference solution at $l_{\rm max\_comp}=8$. Note that the numerical errors are estimated from the acoustic signal that propagates mostly on the very coarse baseline mesh, the impact of the mesh refinement taking only effect inside the vortex structures where the acoustic waves are generated. Consequently, it is difficult to estimate a convergence rate from the overall solution and this explain why $\varepsilon_{p}$ in Table~\ref{tab:summary_2D_L2_pressure} does not follow a second order rate of convergence as in the previous canonical test cases.

Recall that $L$ is the total number of levels of the multi-levels grid set when the hybrid method is employed. As shown in Figure~\ref{fig:results_2D_aeroacoustics_pressure_hybrid2}, solving the fully compressible equations with  $l_{\rm max\_comp}=4$ provides an inaccurate solution for the acoustic pressure. The contribution of $1$ additional layer where the low-Mach-number equations are solved helps to get a pressure field similar to the purely compressible solution computed with $l_{\rm max\_comp}=6$. As reported in Table~\ref{tab:summary_2D_L2_pressure}, simulations with the hybrid method on $L=5$ total levels provide a similar error than the purely compressible approach with $l_{\rm max\_comp}=6$. Furthermore, when the low-Mach-number equations are solved on $2$ additional layers of mesh refinement, i.e.  $L=6$ total levels, the hybrid method recovers the purely compressible solution computed with $l_{\rm max\_comp}=7$. 

\begin{table}
\centering
\renewcommand\arraystretch{1.5}
\begin{tabular}{c || c | c | c | c | c  }
\renewcommand*{\arraystretch}{0.5}
$l_{\rm max\_comp}$ & $L$ & $\varepsilon_{p}$ &  $\Delta t_{\rm LM}$ [s] & $\Delta t_{\rm Comp}$ [s] or $\Delta t_{\rm Hyb}$ [s] & Computational time [s]  \\ \hline\hline 
$1$ &  $\times$ & $\times$ & $\times$ & $2.75 \times 10^{-6}$ &  $13.6$ \\ \hline
$2$ &  $\times$ &  $9.50 \times 10^{-1}$ & $\times$ & $1.37 \times 10^{-6}$ & $54.4$ \\ \hline
$3$ &  $\times$ & $4.53 \times 10^{-1}$ & $\times$ & $6.88 \times 10^{-7}$ &  $240$ \\ \hline
$4$ &  $\times$  & $3.52 \times 10^{-1}$ & $\times$ & $3.44 \times 10^{-7}$ &  $1112$ \\ \hline
$5$ &  $\times$  & $2.45 \times 10^{-1}$ & $\times$ & $1.72 \times 10^{-7}$ &  $4880$ \\ \hline
$6$ &  $\times$  & $1.35 \times 10^{-1}$ & $\times$ & $8.60 \times 10^{-8}$ &  $41080$\\ \hline
$7$ &  $\times$  & $0.57 \times 10^{-1}$ & $\times$ & $4.30 \times 10^{-8}$ &  $303016$ \\ \hline
$4$ &  $5$ & $1.41 \times 10^{-1}$ & $2.7 \times 10^{-6}$ & $3.37 \times 10^{-7}$ & $4936$ \\ \hline
$4$ &  $6$ & $0.69 \times 10^{-1}$ & $1.35 \times 10^{-6}$ & $3.37 \times 10^{-7}$ &  $14880$ \\
\end{tabular}
\caption{Results for the $\mathcal{L}^2$-norm error in the pressure fluctuations $p_1$ ($\varepsilon_{p}$), wall-clock computational time and different time-steps involved in simulations performed with the purely compressible approach and the hybrid method, and for different levels of refinement.}
\label{tab:summary_2D_L2_pressure}
\end{table}

An interesting result here is that the hybrid method is able to recover the purely compressible solution with fewer total levels. This represents a gain in terms of computational burden. Moreover, as $l_{\rm max\_comp} < L$ with the hybrid method, there is also a gain in the time-step. The central columns in  Table~\ref{tab:summary_2D_L2_pressure} present the averaged time-steps $\Delta t_{\rm LM}$ and $\Delta t_{\rm Comp}$ for each simulation performed. Note that when the hybrid method is employed, $\Delta t_{\rm Hyb}$ is reported. The wall-clock CPU time spent for each simulation to reach $t=4 \times 10^{-3}$~s is also reported in the right column. It is interesting to notice that the hybrid method with $l_{\rm max\_comp}=4$ and $1$ additional low-Mach-number level (i.e. $L=5$), the computational time is fairly the same as a purely compressible simulation with $l_{\rm max\_comp}=5$. However the error $\varepsilon_{p}$ corresponds to a purely compressible simulation with $l_{\rm max\_comp}=6$, which means that for a similar solution the hybrid method is about $8.4$ times faster than the purely compressible approach. More interesting, when the simulation is computed with the hybrid method with $l_{\rm max\_comp}=4$ and $2$ additional low-Mach-number levels (i.e. $L=6$), the computational time is about $2.75$ times faster than a purely compressible simulation with $l_{\rm max\_comp}=6$, but as the error $\varepsilon_{p}$ corresponds to a purely compressible simulation with $l_{\rm max\_comp}=7$, the hybrid method is about $7.5$ times faster than the purely compressible approach, which represent a significant gain in the computational time.

\section{Conclusions}
\label{sec:conclusions}

A novel hybrid strategy has been presented in this paper to simulate flows in which the primary features of interest do not rely on high-frequency acoustic effects, but in which long-wavelength acoustics play a nontrivial role and present a computational challenge. Instead of integrating the whole computational domain with the purely compressible equations, which can be prohibitively expensive due to the CFL time step constraint, or with only the low-Mach-number equations, which would remove all acoustic wave propagation, an algorithm has been developed to couple the purely compressible and low-Mach-number equations.   In this new approach, the fully compressible Euler equations are solved on the entire domain, eventually with local refinement, while their low-Mach-number counterparts are solved on specific sub-regions of the domain with higher spatial resolution. The coarser acoustic solution communicates inhomogeneous divergence constraints to the finer low-Mach-number grid, so that the low-Mach-number method retains the long-wavelength acoustics. This strategy fits naturally within the paradigm of block-structured adaptive mesh refinement (AMR) and the present algorithm is developed within the BoxLib framework that provides support for the development of parallel structured-grid AMR applications.

The performance of the hybrid algorithm has been demonstrated on a series of test cases. The temporal and spatial rates of convergence have been investigated with two test cases: first, the propagation of acoustic waves in a uni-dimensional domain; second, the combination of mixed modes composed of the propagation of a circular acoustic wave together with the convection of an entropy spot superimposed to a circular vortex. It has been shown that the acoustic phenomena must be well enough resolved and that solving the low-Mach-number equations on additional levels of mesh refinement helps to get a better solution on other flow phenomena not directly related to the acoustics.

The third test case consists of the simulation of a Kelvin-Helmholtz instability in low-Mach-number mixing layers, which is representative of realistic physical phenomena that can appear in practical flow applications. The initial flow is low-Mach-number and is perturbed so as to generate the formation of vortices that eventually merge together, generating sources of pressure that propagate in the domain. As demonstrated in the present paper, the mixing layer interface requires fine resolution to accurately capture the acoustics, whose long wavelength does not require such a fine resolution. The hybrid method is applied to this problem, and it is demonstrated that the hybrid method is able to provide a very similar solution compared to a fully compressible approach, 
but with fewer levels of refinement and with a significant gain of about two orders of magnitude in time on the global time-step, leading globally to gain of approximately $8$ on the computational time.

Finally, the hybrid method presented in this paper is a first step in the development of a new kind of algorithm to solve problems that feature a large discrepancy in spatial and temporal scales within the same domain. This opens the way to efficient simulations of complex and multi-physics problems such as combustion instabilities in industrial configurations.

\section*{Acknowledgments}

The work here was supported by the U.S. Department of Energy, Office of Science,
Office of Advanced Scientific Computing Research, Applied Mathematics program under
contract number DE-AC02005CH11231. Part of the simulations were performed using resources of the National Energy Research Scientific Computing Center (NERSC), a DOE Office of Science User Facility supported by the Office of Science of the U.S. Department of Energy under Contract No. DE-AC02-05CH11231.





\bibliographystyle{model1-num-names}
\bibliography{bib_manu_2016}







\end{document}